%% file: Preprint-Of_Sullivan_Models,_Massey_products_and_Twisted_Pontrjagin_Products.tex
\documentclass[9pt,lettersize,leqno,fullpage]{article}
\usepackage[usenames,dvipsnames]{color}
\usepackage{amsmath, amssymb, amscd, amsthm, amsfonts, amsbsy, appendix}
\usepackage{bbm, stmaryrd, mathrsfs, bm, linearb, upgreek, skak, textcomp}
\usepackage{mathtools, titlesec, extarrows, tocloft, url, dsfont}
\usepackage[all]{xy} 
\usepackage{graphicx, pgf}

\input{rgb}
\xyoption{web}
\usepackage[T1]{fontenc}

\DeclareMathAlphabet{\mathpzc}{OT1}{pzc}{m}{it} 
\theoremstyle{plain}
\newtheorem{thm}{Theorem}[section]
\newtheorem{cor}[thm]{Corollary}
\newtheorem{lmm}[thm]{Lemma}
\newtheorem{prpn}[thm]{Proposition}
\newtheorem{rem}[thm]{Remark}

\theoremstyle{definition}
\newtheorem{defn}[thm]{Definition}
\newtheorem{eg}[thm]{Example}

\newtheorem{co}[]{}
\newtheorem*{egi}{Example I}
\newtheorem*{egii}{Example II}

\theoremstyle{plain}
\newtheorem*{thma}{Theorem A}
\newtheorem*{thmb}{Theorem B}
\newtheorem*{prpnu}{Proposition a}
\newtheorem*{prpnbu}{Proposition b}

\numberwithin{equation}{section}

\newcommand{\ls}{\hspace*{0.02cm}:}

\include{Macros_Master}

\title{{\LARGE Of Sullivan models, Massey products, and twisted Pontrjagin products}}
\author{Somnath Basu}
\date{}
\begin{document}

\maketitle

\begin{abstract}
Associated to every connected, topological space $X$ there is a Hopf algebra - the Pontrjagin ring of the based loop space of the configuration space of two points in $X$. We prove that this Hopf algebra is not a homotopy invariant of the space. We also exhibit interesting examples of $H$-spaces, which are homotopy equivalent as spaces, which lead to isomorphic rational Hopf algebras or not, depending crucially on the existence of Whitehead products. Moreover, we investigate a (naturally motivated) twisted version of these Pontrjagin rings in the various aforementioned contexts. In all of these examples, Massey products abound and play a key role. 
\end{abstract}
\vspace*{0.7cm}
\tableofcontents
\vspace*{1.5cm}

%=================================================== 1 Introduction ========================================================

\section{Introduction}

\hf\hf Given two connected topological spaces $X$ and $Y$, a few natural questions, if applicable, present themselves\ls
\bgd
\textup{{\it Are the two spaces diffeomorphic\,? homeomorphic\,? homotopy equivalent\,?}}
\edd
{\it Lens spaces} provide the first non-trivial examples of spaces that are homotopy equivalent but not homeomorphic. Milnor's celebrated exotic spheres provided the first examples of manifolds that are homeomorphic but not diffeomorphic! In order to detect the difference between two spaces, one first checks the usual suspects like homotopy groups, and (co)homology groups. If these are isomorphic then the structure of the cohomology ring, arising from the cup product, may resolve the stalemate. For instance, $S^2\times S^2$ and $\mathbb{CP}^2\#\overline{\mathbb{CP}}^2$ have identical homotopy groups and homology groups. Their integral cohomology rings are, however, not isomorphic. But now imagine a scenario where cohomology rings aren't sufficient\ls\hspace*{0.05cm}Consider the complement, in $\R^3$, of the links given below.
\begin{figure}[h]
\begin{center}
\includegraphics[scale=0.5]{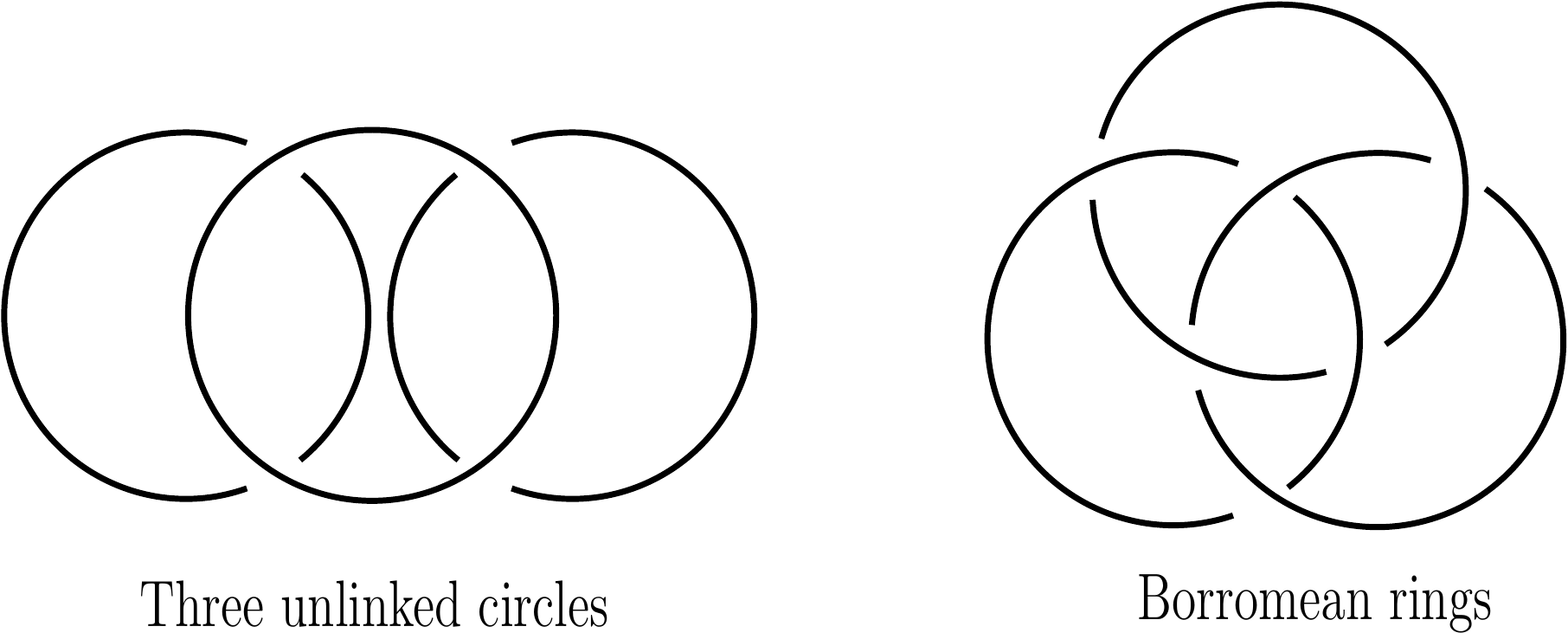}\vspace{-0.2cm}
\caption{Two different ways three circles interact in $\R^3$.}
\end{center}\vspace*{-0.2cm}
\end{figure}
Notice that although any two circles in the Borromean rings\footnote{The name Borromean rings comes from their use by the Borromeo family in Italy. The existence and its use, however, goes back to the 2nd century in Buddhist cultures.} are unlinked, the three taken together are not! One should contrast this with the three totally unlinked circles on the left. Intuitively, one {\it knows} that the associated  complementary spaces are not the same. One can show that these spaces cannot be distinguished from each other by usual invariants, and their cohomology rings are isomorphic. However, there are secondary invariants called {\it Massey products} which makes our intuition precise; it distinguishes the two spaces under consideration by capturing this {\it triple linking}, in this instance. For a lovely exposition of this and more, we highly recommend Massey's original classic \cite{Mas69}. \\
\hf\hf Yet another approach, which has strong correlation with existence of Massey products, is to study the configuration space $F_n(M)$ of $n$ points of a given space $M$, which is the {\it state space} of $n$ particles moving around in $M$. In particular, $F_2(M):=(M\times M)\setminus M$ is the configuration space of two points. Configuration spaces are extremely rich in structure and have been extensively studied. It has strong connections with homotopy theory and robotics, among other topics. As a particular instance, it was unknown for a long time, but largely believed to be true, if the following holds\ls
\bgd
\textit{If $M$ and $N$ are closed, homotopy equivalent manifolds then $F_2(M)$ is homotopy equivalent to $F_2(N)$.}
\edd
It turns out to be false! One could imagine, albeit naively, that the $2$-point configuration space probes the manifold, captures the neighbourhood of the diagonal inside $M\times M$, and thereby produce invariants subtle enough to be not detected by homotopy equivalences. It was shown by Levitt \cite{Lev95} that the based loop space of the configuration space preserves homotopy equivalences. In 2005 Longoni and Salvatore \cite{LS05} exhibited homotopy equivalent lens spaces (and necessarily non-homeomorphic for the purpose) which have non-homotopic configuration spaces. In fact, they use Massey product on the universal cover of configuration spaces to distinguish the two. \\
\hf\hf Let us now consider another example\ls\hspace*{0.05cm}distinguish the unitary group $SU(3)$ from $S^3\times S^5$. The homology groups are identical and the cohomology ring structures match. What about homotopy groups? If we ignore torsion (or pretend not to have prior knowledge of higher homotopy groups of spheres) then it can be shown that there is a map which induces an isomorphism barring torsion elements. We know that $SU(3)$ is a Lie group while $S^3\times S^5$ isn't one. Moreover, there is a canonical multiplicative map, arising from the group multiplication, on the homology of $SU(3)$ which is absent for $S^3\times S^5$. This leads us to spaces that have more structure.\\
\hf\hf We now assume that the spaces $X$ and $Y$ come equipped with a {\it group-like} multiplication, i.e., there is a map $m_X:X\times X\to X$ that is associative only up to homotopy. Such spaces are called $H$-spaces, the name coined by Serre in honour of Heinz Hopf. Notice that this multiplicative structure induces a product on homology of the space. We call this the {\it Pontrjagin product} associated to $X$, and the homology with this product structure is called the {\it Pontrjagin ring} of $X$. The canonical example of such spaces via homotopy theory is that of a based loop space of a (topological) space, where the Pontrjagin product is given by concatenation of loops.\\
\hf\hf Mathematical sensibility dictates that the natural question to ask here is not if $X$ and $Y$ are homotopic (which they very well may be), but rather if the spaces are homotopic as $H$-spaces. In particular, is there a homotopy equivalence $f:X\to Y$ which preserves the Pontrjagin product? We thus make a full circle and come back to our original question in the context of $H$-spaces. Consider the question\ls
\bgd
\textup{{\it Are the based loop spaces of the configuration spaces of homotopy equivalent manifolds equivalent as $H$-spaces?}}
\edd
This is a natural question that was not addressed by Levitt \cite{Lev95}. We explore this question, and prove, among other things, that the answer to the question above is negative! More generally, we focus on structures that play a key role in distinguishing two $H$-spaces $X$ and $Y$ which are homotopy equivalent as spaces. Put simply, we want to analyze examples where the Pontrjagin product associated with $X$ is not isomorphic to that of $Y$. Along the way, we exhibit a strong interplay between the Pontrjagin product and Massey products. In fact, this connection is subtle and actually stems from the existence of {\it Whitehead products}. Our main language to process all these structures is that of {\it Sullivan models} and {\it rational homotopy theory}. \\
\hf\hf To further our scheme, we digress into a discussion on fibrations and their associated based loop fibrations. Let $F\into E\to B$ be a {\it fibration}\footnote{By a fibration we mean a {\it Serre fibration}, i.e., one which has the usual homotopy lifting property.} where $F,E$, and $B$ are connected topological spaces. Applying the based loop functor to the original fibration yields a fibration 
\bgd
\Omega F\into \Omega E\to \Omega B,
\edd
where $\Omega X$ is the space of (continuous) loops in $X$ based at $x_0\in X$. Notice that $\Omega X$ is an $H$-space with concatenation of loops, denoted by $\ast$. If we have a section $s:B\to E$ then there is a natural map 
\bgd
\mathpzc{s}:\Omega B\times \Omega F\longrightarrow \Omega E,\,\,(\gamma_1,\gamma_2)\mapsto s(\gamma_1)\ast\gamma_2.
\edd
The map above need not be a map of $H$-spaces as can be seen from the following\ls
\begin{eg}
Recall that $K(G,1)$, for a discrete group $G$, is a connected topological space which has precisely one non-zero homotopy group and
\bgd
\pi_1(K(G,1))\cong G.
\edd
This space is uniquely determined up to homotopy equivalence and is commonly referred to as the {\it classifying space} of $G$. It is also an example of an {\it Eilenberg-Maclane space}. Start with any short exact sequence of (discrete) groups $0\to A\to B\to C\to 0$, equipped with a section $s:C\to B$. Using $B\to C$ form a fibration $p:K(B,1)\to K(C,1)$ of {\it Eilenberg-Maclane spaces} with the homotopy fibre $K(A,1)$. Of course, if $B\not\cong A\times C$ then the map at the level of fundamental groups does not split. The section induces a map
\bgd
\mathpzc{s}:\Omega K(A,1)\times \Omega K(C,1)\longrightarrow \Omega K(B,1)
\edd
which is a homotopy equivalence as a map of spaces. However, if we identify $\Omega K(G,1)$ with $G$, then $s$ is a not a homotopy equivalence in the category of (discrete) groups, i.e., the map is {\it not} a group isomorphism! 
\end{eg}
When $F,E$, and $B$ are simply connected $\mathpzc{s}$ would induce a homotopy equivalence. However, the nagging truth persists\ls\hspace*{0.05cm}This may not be a map of $H$-spaces. The Pontrjagin product may still be powerful enough, in a class of examples which we outline subsequently, to distinguish between the $H$-spaces $\Omega B\times \Omega F$ and $\Omega E$. \\[0.2cm]

\bgc
{\bf Connections with String Topology}
\edc
\hf\hf We would like to emphasize the link between the contents of this article and {\it string topology}, the study of differential and algebraic topology of free loop spaces of manifolds. We refer the interested reader to \cite{CS99}. There are a couple of different approaches to this subject and quite a few relevant literature on this topic. In any case, a lot of algebraic structures came forth in this study and it was asked if some of these distinguish between homeomorphic manifolds which are not diffeomorphic. One can even ask if these structures help differentiate between homotopy equivalent but non-homeomorphic manifolds. Most of the usual algebraic invariants arising from string topology were shown to be homotopy invariants.\\
\hf\hf There is, however, a geometrically altered version of string topology that goes by the name of {\it transversal string topology\footnote{As Whitehead had wisely said, "`Transversal' is a noun. The adjective is `transverse'." We understand this but we'll stick to transversal for now.}} and was studied in the author's Ph.D. thesis. The subject of transversal string topology and its underpinnings will be taken up elsewhere. A rather interesting property of transversal string topology is that it is {\it not} a homotopy invariant. A crucial input that goes in proving this is part (ii) of {\bf Theorem B}. The twisted versions of the main results (part (ii) of {\bf Theorem A} and {\bf Theorem B}) are part of a bigger picture which naturally lives in the world of string topology. Finally, the results of {\bf Theorem A} are also inspired by string topology. \\[0.2cm]

\bgc
{\bf Main Results}
\edc
\hf\hf We briefly present below the key results ({\bf Theorem A}, {\bf Theorem B}, {\bf Proposition a}, and {\bf Proposition b}) of the paper in the context of the two examples we shall study, fitting ideally in the framework of $H$-spaces. It is intentionally stated in less symbolic terms for the ease of the reader. The more technical and symbolic versions are referenced for those who wish to jump ahead and read it instead. 
\begin{egi}\label{freeloop}
Consider the {\it free loop space} $LM:=C^0(S^1,M)$ associated to a $2$-connected\footnote{A space is called $2$-connected if it is connected and the first two homotopy groups vanish.} manifold $M$. The natural evaluation map (evaluating a loop at $1\in S^1)$ defines a fibration $\textup{ev}:LM\to M$ with fibre the based loop space $\Omega M:=C^0((S^1,1),(M,x_0))$. This fibration has a natural section $s:M\to LM$ by sending a point to the constant loop at that point. Moreover, the $2$-connectedness of $M$ ensures that all the spaces in the fibration are simply connected. In general, $\Omega M\times M$ is very far from being homotopic to $LM$. However, the induced map
\bgd
\mathpzc{s}:\Omega M\times \Omega (\Omega M)\longrightarrow \Omega (LM)
\edd
is a homotopy equivalence. Are these different as $H$-spaces? We have an affirmative answer when $M$ is an even dimensional sphere (cf. Theorem \ref{LLSS})\ls
\begin{thma}
Let $S^{2l}$ denote the sphere of dimension $2l$ for $l\geq 1$. Consider the two spaces :\\
\hf $\Omega(LS^{2l})$ - the based loop space of the free loop space of $2l$-sphere, and\\
\hf $\Omega(S^{2l}\times\Omega S^{2l})$ - the based loop space of the product of $S^{2l}$ and the based loop space of $S^{2l}$. 
\begin{itemize}
\item[\textup{(i)}] 
The two spaces above are homotopy equivalent but not homotopic as $H$-spaces.
\item[\textup{(ii)}] Consider the spaces $\Omega(LS^{2l})$ and $\Omega(S^{2l}\times\Omega S^{2l})$ but with a twisted multiplication (cf. Definition \ref{twmul}) on both. These are not homotopic as twisted $H$-spaces.
\end{itemize}
\end{thma}
For the interested reader, we twist the original multiplication by a generator of the group $\pi_{2l-2}$, which is $\Z$ for both the spaces above. The existence of Massey products for $LM$ (cf. Theorem \ref{MasseyS2l}) plays a crucial role. However, the subtlety lies in the fact that it is the quadratic nature of the differential of the minimal model of $LM$ that leads to non-trivial {\it Whitehead products}, interesting Massey products, and to seemingly incongruous Pontrjagin products. As an illustration of this, we prove that Theorem A (cf. Theorem\ref{LLSS}) fails to hold when $M=\mathbb{CP}^n$ and $n\geq 2$ (cf. Propostion \ref{OLCP})\ls
\begin{prpnu}
There is a natural map 
\bgd
\Phi:\Omega (\mathbb{CP}^n\times\Omega \mathbb{CP}^n)\longrightarrow\Omega (L\mathbb{CP}^n),
\edd
and it is a rational homotopy equivalence of $H$-spaces if $n\geq 2$.
\end{prpnu}
\end{egi}
\begin{egii}\label{lens}
Let $M$ be any compact, oriented manifold which has a nowhere vanishing vector field $V$. Let $F_2(M)$ denote the $2$-point configuration space of $M$. This space maps to $M$ (by projecting on to the first coordinate) with fibre $M-\{x_0\}$. It is known that this is a locally trivial fibre bundle. The existence of $V$, using the exponential map, implies the existance of a section $s(x)=(x,\textup{exp}(\ep V(x)))$ for $\ep>0$ small enough. Then we have a map
\bgd
\mathpzc{s}:\Omega M\times \Omega (M-\{x_0\})\longrightarrow \Omega(F_2(M)),
\edd
which is a homotopy equivalence. However, as the example of lens spaces (cf. \S \ref{Conf2}, \S \ref{Conf22}) $L_{7,1}$ and $L_{7,2}$ shows, this map is not a map of $H$-spaces. This follows, for instance, from the fact that $\Omega F_2(L_{7,1})$ and $\Omega F_2(L_{7,2})$ have non-isomorphic Pontrjagin products (cf. Corollary \ref{ueac1}, Theorem \ref{L72}, and Theorem \ref{L72tw})\ls
\begin{thmb}\label{thmB}
Let $\overline{X}$ denote the universal cover of a topological space $X$. Then the following holds\ls
\begin{itemize}
\item[\textup{(i)}] The Pontrjagin rings associated to $\Omega \overline{F_2(L_{7,1})}$ and $\Omega\overline{F_2(L_{7,2})}$ are non-isomorphic.
\item[\textup{(ii)}] The twisted Pontrjagin rings associated to $\Omega \overline{F_2(L_{7,1})}$ and $\Omega\overline{F_2(L_{7,2})}$ are non-isomorphic.
\end{itemize}
\end{thmb}
We prove this by exhibiting central elements associated to $L_{7,1}$ which are absent for $L_{7,2}$. Again, for the interested reader, the twisting is by a natural element in the fundamental group that corresponds to the normal sphere of a point (on the diagonal) inside $L_{7,j}\times L_{7,j}$. On the other hand, the theorem above holds true (in both the twisted and the untwisted version) for based loop spaces of higher configuration spaces (cf. Proposition \ref{lensn} and the concluding remarks)\ls
\begin{prpnbu}
For $n>2$ the based loop spaces $\Omega \overline{F_n(L_{7,1})}$ and $\Omega \overline{F_n(L_{7,2})}$ are not homotopy equivalent as $H$-spaces or as twisted $H$-spaces.
\end{prpnbu}
\end{egii}
\hf\hf Having announced the main results, we find it necessary to mention that the proofs illustrate the propagation of Massey products into Pontrjagin products. We also briefly explain why we have included a twisted version of the same. The proverbial twist in the tale, comes from studying certain twisted (cf. Definition \ref{twmul}) Pontrjagin products, and discovering that the Massey products still persist in this setting! The proof is essentially similar to the untwisted case. The reason for the inclusion of the twisted version of the results is that these appear naturally in the setting of {\it transversal string topology}. This topic will be subsequently taken up elsewhere in forthcoming papers and was part of the author's dissertation. The twisted Pontrjagin product appears as a by-product of applying the transversal string topology machinery and is the final ingredient towards proving that transversal string topology is {\it not} a homotopy invariant. However, we felt that the results outlined here are of considerable independent interest.\\[0.5cm]
\hf\hf {\bf Outline of the paper}\ls\hspace*{0.05cm}In \S \ref{Mplay} we give a brief introduction to the triad of mathematical ideas - Sullivan models, Massey products, and Pontrjagin products. In a way this paper can be thought of as a word in these three generators. In \S \ref{MassConf} we analyze Example II by introducing the relevant background material and translating everything in terms of minimal models, and processing the subsequent computations. In \S \ref{MassLM} we analyze Example I by reviewing the minimal model for free loop spaces. We show that Massey products exist, and subsequently gives rise to interesting Pontrjagin products, or not, depending on the non-triviality of Whitehead products. It should be kept in mind that almost all computations of the Pontrjagin product relies crucially on the famous result (cf. Thereom \ref{MM}) of Milnor-Moore \cite{MM65} which allows us to do the computations in the {\it universal enveloping algebra} of the Lie algebra generated by the appropriate homotopy groups.\\[0.5cm]
\hf\hf {\bf Acknowledgements\ls} This paper would not have been possible without the help and direction of Professor Dennis Sullivan. \\[0.5cm]

%=============================== 1.1 Minimal models, Milnor-Moore, Massey product ==========================================

\subsection{The main players}\label{Mplay}

\hf\hf In the course of this article, three main mathematical ideas will be playing key roles. It's the interplay between the trio that is of most interest to us. These ideas have been named already in the title of this paper. We shall begin by a brief review of these concepts. The experienced reader may choose to skip this section and refer to this as and when needed. We thought it would be useful for a general reader to be familiar with these before they delve into the technical details of the paper, as these concepts make quite frequent appearance throughout.\\

%================================================= Minimal Models ===============================================

\begin{co}\label{MinMod}{\bf Sullivan Models}\\  
 
\hf\hf Spaces can be modelled algebraically in the sense of rational homotopy theory. Fibrations of spaces can be treated similarly - we model the base space, and then model the twisting of the fibration appropriately. These algebraic models, usually called {\it minimal models} and often called {\it Sullivan models}, was originally introduced by D. Sullivan and is built out of piecewise linear differential forms. The classic \cite{SulIHES} by Sullivan provides a coherent and unifying approach in describing various naturally occuring objects in terms of minimal models. We briefly review the existence of such a model for any simply connected manifold. The standard (and almost comprehensive) reference for rational homotopy theory is the book \cite{FHT01}. The theory of minimal models is an extremely powerful tool, and has had a strong influence on the significance of doing homotopy theory over $\Q$, ever since its introduction. It has also been used to prove beautiful theorems, \cite{DGMS75} and \cite{VPS76} to name a few. In this discussion all subsequent models are rational models, i.e., we will work with rational coefficients unless otherwise stated.\\
\hf\hf Let $\Lambda V$ be the free graded commutative algebra generated by a graded $\R$-vector space $V=\oplus_{i\geq 0} V_i$. When $\Lambda V$ is equipped with a differential $d$, we call $(\Lambda V,d)$ a {\it minimal algebra} if there is a set of generators $\{a_i\}$, indexed by a well ordered set $I$, such that $|a_i|\leq |a_j|$ if $i<j$ and $da_i$ can be written in terms of $a_j$'s for $j<i$. In other words, $da_i$ doesn't have a {\it linear term}. A minimal algebra $(\Lambda V,d)$ is said to be a {\it minimal model} for a connected manifold $M$ if there is a quasi-isomorphism\footnote{In general, when one speaks of a quasi-isomorphism from $F$ to $W$ one always assumes it is given by a {\it zig-zag}
\bgd
F\to W_1\leftarrow F_1\to\cdots\to F_n\to W
\edd
where each arrow induce isomorphism in (co)homology. In the case of minimal models, if $F_i,W_i$'s and $F$ are all minimal then all the arrows are invertible and hence we make no mention of zig-zag when dealing with quasi-isomorphisms.} of algebras $\varphi:(\Lambda V,d)\to (\Omega_{dR}M,d)$, where $\Omega_{dR}M$ is the space of de Rham forms. These discussions can also be carried out with rational coefficients and replacing $\Omega_{dR}M$ with 
\bgd
C^{\lp}(M;\Q):=\oplus_{i\geq 0}C^i(M;\Q),
\edd
the cochain complex of $M$. It is possible, some times, to have a quasi-isomorphism between the minimal model and $(H^{\lp}(M;\Q),0)$, i.e., the graded cohomology ring with the zero differential. In such a case, $M$ is called {\it formal}. Henceforth, we shall assume that $V$, as above, is a rational vector space. \\
\hf\hf It is known that minimal models exist for any connected topological space, and is unique up to isomorphism (\cite{FHT01}, Theorem 14.12). From Sullivan's work on rational homotopy theory, explicit models can be given for any {\it nilpotent space} $X$, i.e., $\pi_1(X)$ is nilpotent, and its action on the higher homotopy groups is nilpotent. We shall work with simply connected spaces, where the dual of the rational homotopy groups $\pi_i^{\Q}(X):=\pi_i(X)\otimes\Q$ can be taken to be $V_i$. Moreover, the quadratic part of $d$ is dual to the Whitehead product (cf. Definition \ref{Wh}) on $\pi^{\Q}_{\lp}(X):=\oplus_{i\geq 2}\pi_i^{\Q}(X)$. For a detailed discussion on this and more, we refer the reader to \cite{FHT01}. Some examples of such models are in order.
\begin{eg}\label{sph}{\bf (Spheres)}\\[0.2cm]
For $S^{2k+1}$ we take the free (graded commutative) algebra $(\Lambda(x),d\equiv 0)$ generated by a generator $x$ in degree $2k+1$. The map 
\bgd
\varphi:(\Lambda(x),d)\stackrel{\simeq}{\longrightarrow} H^{\lp}(S^{2k+1};\Q)
\edd
is given by mapping $x$ to the volume form. Similar models can be built for $S^{2k}$ by taking 
\bgd
(\Lambda(x,y),d),\,\,dx=0,dy=x^2,|x|=2k.
\edd
Since the underlying vector spaces of the minimal model can also be built out of (rational) homotopy groups, one recovers part of the celebrated result of Serre on homotopy groups of spheres\ls \\
\hf (i) $\pi_i^{\Q}(S^{2k+1})=\Q$ if $i=0,2k+1$ and zero otherwise,\\
\hf (ii) $\pi_i^{\Q}(S^{2k})=\Q$ if $i=0,2k,4k-1$ and zero otherwise.
\end{eg}
\begin{eg}{\bf (Complex Projective Spaces)}\\[0.2cm]
Let $\mathbb{CP}^n$ be the complex projective space of (complex) dimension $n$. The cohomology ring 
\bgd
H^{\lp}(\mathbb{CP}^n;\Q)\cong\Q[\alpha]/(\alpha^{n+1})
\edd
is a truncated algebra generated by an element of degree two. A minimal model is given by $(\Lambda(x,y),d)$, where
\bgd
dx=0,\,\,dy=x^{n+1},\,\,|x|=2.
\edd
The map 
\bgd
\varphi:(\Lambda(x,y),d)\stackrel{\simeq}{\longrightarrow} (H^{\lp}(\mathbb{CP}^n;\Q),0),
\edd
sending $x$ to $\alpha$, and $y$ to $0$, is a quasi-isomorphism. 
\end{eg}
\end{co}

%================================================= Massey Products ==============================================

\begin{co}\label{MassProd}{\bf Massey Products}\\   

\hf\hf Let us revisit the problem of deciding if two given topological spaces $X$ and $Y$ are homotopy equivalent. If all the natural invariants turn out to be useless at this task and even the cohomology rings are isomorphic, one may consider secondary\footnote{The terminology derives its name from the fact that the cohomology ring of a space is a primary invariant.} invariants called {\it Massey products} (cf. \S \ref{MassProd}) which often distinguish $X$ from $Y$. In particular, it is very useful when the cup product structure is zero while this invariant is non-zero. In Massey's original paper \cite{Mas69}, where Massey products are first defined, Massey calls it {\it higher order linking numbers}. The example of Borromean rings, as explained in the introduction, provides the answer to this nomenclature. The classical example, as outlined by Massey himself, is far too pretty to go unread by a modern mathematician, and the author urges the reader, again, to read the classic \cite{Mas69}. Extremely geometric in its nature, Massey uses his invariant to detect higher-order linking in Borromean rings, i.e., any two taken together are unlinked but the three taken together cannot be unlinked! A homological picture of the same is provided below\ls
\begin{figure}[h]
\begin{center}
\includegraphics[scale=0.5]{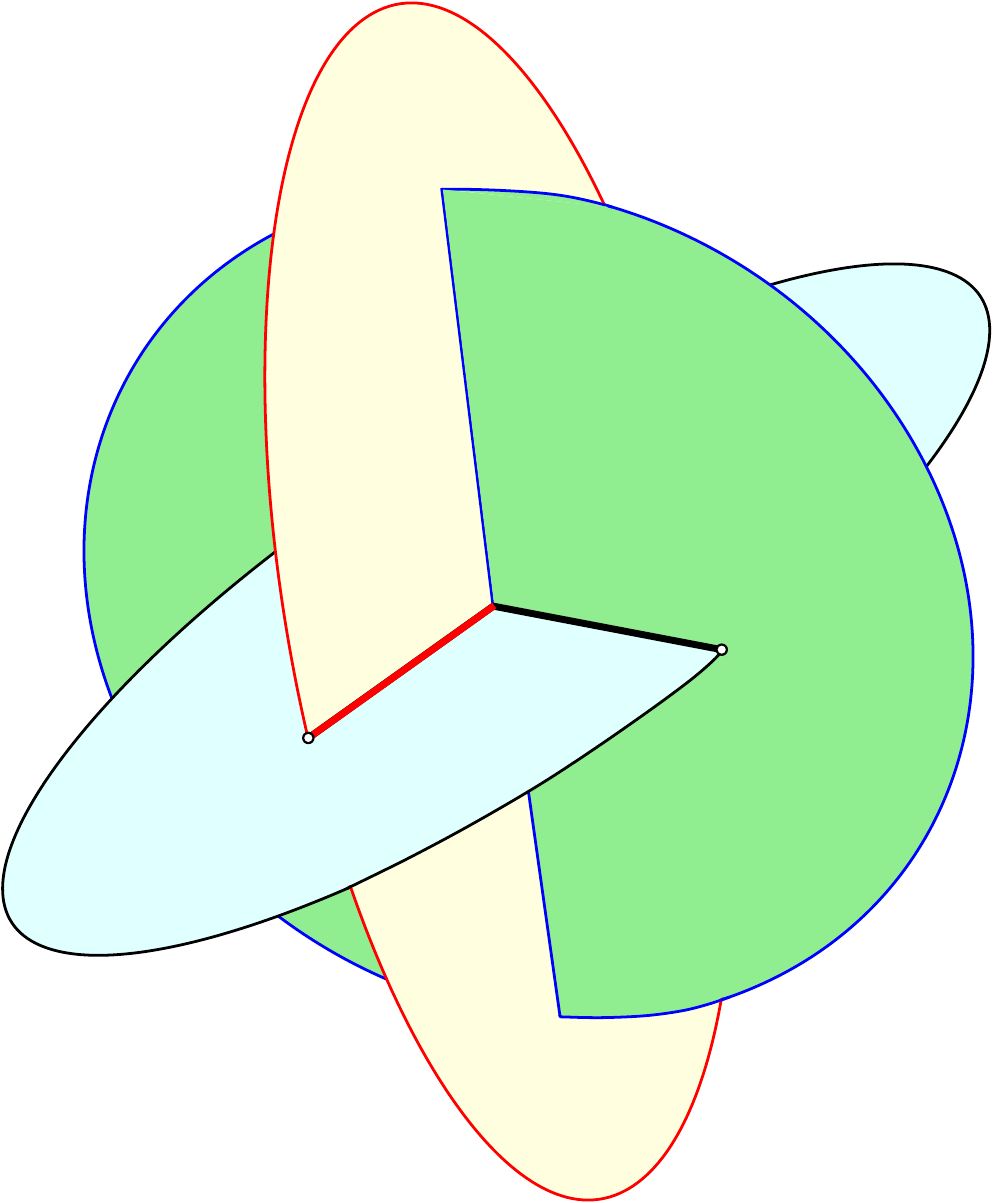}\vspace{-0.2cm}
\caption{The triple linking captured by the transversal intersection of three disks.}
\end{center}\vspace*{-0.2cm}
\end{figure}
In fact, a homological version of Massey product is presented in \cite{Mas69} and we're going to use it later. \\
\begin{defn}{\bf (Triple Massey products)}\\
Let $(A,d)$ be a differential graded algebra with $H^{\lp}(A)$, its graded cohomology ring. For pure elements $u,v,w\in H^{\lp}(A)$, choose representatives $U,V,W\in A$. If $UV=dS$ and $VW=dT$ then define the {\it Massey triple product} of $u,v,w$ to be
\bge\label{Massey}
\left\langle u,v,w\right\rangle:=(-1)^{|v|}[SW-(-1)^{|u|}UT]\in H^{\lp}(A)/(u,w),
\ede
where $(u,w)$ is the ideal generated by $u$ and $w$.
\end{defn}
Notice that although $\left\langle u,v,w\right\rangle$ is defined, $\left\langle v,u,w\right\rangle$ may not be. However, there is a graded commutativity\ls
\bgd
\left\langle u,v,w\right\rangle =-(-1)^{|u||v|+|v||w|+|w||u|}\left\langle w,v,u\right\rangle.
\edd
\hf\hf For a (connected) topological space $X$ set $(A,d):=(C^{\lp}(X),\delta)$, the cochain complex of $X$. There are notions of Massey products of higher orders and of order two. The second order Massey product is the usual cup product which is a primary invariant of $X$. The higher order Massey products, up to indeterminacy, are secondary invariants of the space. Manifolds on which all Massey products vanish are called {\it formal} manifolds. The rational homotopy type of such manifolds is determined by its rational cohomology ring. This is true, for instance, in the case of compact, K\"ahler manifolds \cite{DGMS75}. Formality is a key property that provides a connection between the homotopy type of the manifold with other geometric structures that the manifold admit. It is a property that is preserved under products, joins and wedges of spaces. Finally, there are higher order Massey products which detect higher order linking. Among other uses, Massey products appear naturally in {\it twisted $K$-theory}. 
\end{co}

%============================================= Pontrjagin Products ==============================================

\begin{co}\label{PontProd}{\bf Pontrjagin Products}\\   

\hf\hf Let $(X,x_0)$ be a (based) connected topological space and consider $\Omega X$, the {\it based loop space} of $X$. Since we can compose based loops by concatenation, $\Omega X$ is an $H$-space. In any case, we have a map
\bgd
\mathpzc{m}:\Omega X \times\Omega X\longrightarrow\Omega X
\edd
which induces a map on homology $H_{\lp}(\Omega X;\Q):=\oplus_{i\geq 0}H_i(\Omega X;\Q)$ and can be written as
\bge\label{Pont}
\times : H_{\lp}(\Omega X;\Q)\otimes H_{\lp}(\Omega X;\Q)\longrightarrow H_{\lp}(\Omega X;\Q).
\ede
This is called the {\it Pontrjagin product} on the based loop space of $X$, a notion introduced in more generality by Pontrjagin in 1939. We also have the map
\bgd
\Delta : H_{\lp}(\Omega X;\Q)\longrightarrow H_{\lp}(\Omega X;\Q)\otimes H_{\lp}(\Omega X;\Q),
\edd
induced by the diagonal map $\Omega X\to \Omega X\times \Omega X$. These two structure maps turn $H_{\lp}(\Omega X;\Q)$ into a {\it Hopf algebra}. \\
\hf\hf For any topological space $X$, J. H. C. Whitehead had defined a product, now called the {\it Whitehead product}, on the homotopy groups 
\bge\label{Wh}
\textup{Wh}:\pi_i (X)\otimes \pi_j(X)\longrightarrow \pi_{i+j-1}(X),\,\,\forall\,\, i,j\geq 0.
\ede
In fact, Massey and Uehara in 1957 used Massey product to show that the Jacobi identity holds for the Whitehead product. If $X$ is simply connected then the vector space $\pi^{\Q}_{\lp}(X):=\oplus_{i\geq 0} \pi_i^{\Q}(X)$ turns into a graded Lie algebra with $\textup{Wh}$ as the Lie bracket of degree $-1$. Moreover, in such a case, $\Omega X$ is connected, and we have natural isomorphisms $\pi_{i+1}(X)\cong \pi_{i}(\Omega X)$.  If we transfer the Lie algebra structure via these natural maps then we get a Lie algebra structure on $\pi_{\lp}^{\Q}(\Omega X)$, called the {\it Samelson product}. \\
\hf\hf Historically, Samelson had defined this product on the homotopy groups of a topological group and had shown it was bilinear. In a 
topological group $G$ the correspondence $(g,h)\to ghg^{-1}h^{-1}$ defines a map
\bgd
\mathpzc{c}:G\wedge G\longrightarrow G,
\edd
where $G\wedge G$ stands for the {\it smash product} of $G$ with itself, i.e., it is the identification space $G\times G/(G\times\{e\}\cup \{e\}\times G)$. This map $\mathpzc{c}$ induces a pairing of $\pi_i(G)$ with $\pi_j(G)$ to $\pi_{i+j}(G)$ in a natural way. The commutator $\mathpzc{c}$ thus induces a ring structure on $\pi_\lp(G):=\oplus_{i\geq 0}\pi_i(G)$. This is called the {\it Samelson product}. If use the Hurewicz map from the homotopy groups to the homology groups then the Samelson product transfers to the Pontrjagin product.\\
\hf\hf With that brief detour in place, we come back to $\Omega X$. This space, however, is not a topological group but a topological monoid where inverses exist up to homotopy. It's still good enough as one can replace $\Omega X$, due to a construction of Milnor, and independently by Kan, by a topological group. We shall, by abuse of notation, denote both, the Whitehead product and the Samelson product, by $\textup{Wh}$. \\
\hf\hf The observant reader would notice that given a simply connected space $X$, there are two competing structures that naturally present themselves. The first one is a Hopf algebra given by $(H_{\lp}(\Omega X;\Q),\times)$, and the second is the Lie algebra $(\pi_{\lp}^{\Q}(\Omega X),\textup{Wh})$. How are the two structures related? There is a Hopf algebra naturally associated to any Lie algebra, viz., the {\it universal enveloping algebra} of the Lie algebra. Let us recall the following\ls
\begin{defn}{\bf (Universal Enveloping Algebra)}\hf Given a (graded) Lie algebra $(L, [\,,\,])$, the {\it universal enveloping algebra} is an algebra 
$U(L)$ equipped with a morphism of Lie algebras $\iota_L:L\to U(L)$ such that any morphism of Lie algebras $\varphi:L\to A$ factors through $U(L)$, i.e., there exists a morphism of Lie algebras $\widetilde{\varphi}:U(L)\to A$ such that $\widetilde{\varphi}\circ \iota_L =\varphi$. 
\end{defn}
It can be checked that $U(L)$ is unique up to isomorphism. Therefore, an equivalent way of defining $U(L)$ is by providing an explicit presentation.
\begin{defn}\label{UEA2}{\bf (Universal Enveloping Algebra)}\hf Given a (graded) Lie algebra $(L, [\,,\,])$ where the Lie bracket has degree $d$, a presentation of the {\it universal enveloping algebra} is given by $T(L)/I$, where $T(L)$ is the free tensor algebra generated by $L$ and $I$ is the ideal generated by relations of the form
\bge
x\otimes y-(-1)^{(|x|+d)(|y|+d)}y\otimes x-[x,y].
\ede
\end{defn}
\hf\hf The final step would be to connect the homotopy groups of $\Omega X$ to the homology groups of $\Omega X$ with the hope 
that the naturality of this connection generates the isomorphism we are seeking between the two structures. Recall that the {\it Hurewicz map}
\bgd
\mathcal{H}:\pi_i(\Omega X)\longrightarrow H_i(\Omega X)
\edd
provides a natural pathway, and is, in fact, a morphism of Lie algebras. We are now ready to state the beautiful structural result of Milnor \& Moore \cite{MM65}\ls
\begin{thm}\label{MM}{\bf (Milnor-Moore)}\\
If $G$ is a pathwise connected homotopy associative $H$-space with unit, and $\lambda:\pi_\ast(G;K)\to H_\ast(G;K)$ is the Hurewicz morphism of Lie algebras, then the induced morphism $\widetilde{\lambda}:U(\pi_\ast(G;K))\to H_\ast(G;K)$ is an isomorphism of Hopf algebras.
\end{thm}
In the theorem above, $K$ is a field of characteristic zero. As a consequence, for a simply connected space $X$ the conditions of the theorem is obtained for $\Omega X$\ls 
\bge\label{MilMo}
\mathcal{H}:U(\pi_{\lp}^\Q(\Omega X))\stackrel{\cong}{\longrightarrow} \left(H_{\lp}(\Omega X;\Q),\times\right).
\ede
\end{co}

\newpage

%============================================ 2 Massey via configuration spaces ============================================

\section{Massey and Pontrjagin products via configuration spaces}\label{MassConf}

\hf\hf We shall briefly review lens spaces and their classification. We then proceed to study configuration spaces of lens spaces and sketch the method used in \cite{LS05} to prove that configuration spaces doesn't preserve homotopy invariance. We recast their result in the context of minimal models and thereby compute the Pontrjagin product of the based loop space of these configuration spaces. To be exact, we calculate the Pontrjagin ring of $\Omega\overline{X}$ where $\overline{X}$ is the universal cover of $X$, and $X$ is the appropriate configuration space. This structure induces, and is induced by, the Pontrjagin ring structure of $\Omega X$ as the based loop space of $X$ is a union of (homeomorphic) components, indexed by $\pi_1(X)$, where each component is $\Omega\overline{X}$. The Pontrjagin ring of $\Omega X$ is the direct sum of $\pi_1(X)$-copies of the Pontrjagin ring of $\Omega\overline{X}$ and $\pi_1(X)$ acts by on the summands via left translation on the indices labelled by $\pi_1(X)$.

%================================ Section 2.1 - Configuration spaces of two points =========================================

\subsection{Configuration spaces of two points I}\label{Conf2}

\hf\hf The configuration space (of $n$ points) of any space is interesting in its own right. It is the natural domain for a lot of configuration problems as it is the state space of particles moving in the ambient space. It's a well studied object and shows up in numerous places inside mathematics and physics. It's even more intriguing in light of the fact that it contains non-homotopy information, i.e., there are homotopy equivalent spaces $X$ and $Y$ such that the associated configuration spaces are not. The example we'll look at involves lens spaces and is due to Longoni and Salvatore \cite{LS05}. We shall briefly review it and present a minimal model approach to study such configuration spaces. \\
\hf\hf We shall think of $S^3=\{(z_1,z_2)\in \C^2\,|\,|z_1|^2+|z_2|^2=1\}$ as the unit quarternions. We denote by $L_{p,q}$ the {\it lens space} arising as the quotient of $S^3$ by the $p$th roots of unity as follows\ls
\bgd
e^{i2\pi k/p}(z_1,z_2)=(e^{i2\pi k/p}z_1,e^{i2\pi kq/p}z_2).
\edd
Since $S^3$ is the universal cover of $L_{p,q}$, $\pi_1(L_{p,q})={\Z}_p$. The following is well known.
\begin{thm}{\bf (Reidemeister, Brody)}\\
Consider the lens spaces $L_{p,q_1}$ and $L_{p,q_2}$.\\
\hf \textup{(i)} They are homotopy equivalent if and only if $q_1q_2\equiv\pm n^2\,\textup{mod}\,\,p$.\\
\hf \textup{(ii)} They are homeomorphic if and only if $q_1\equiv\pm q_2^{\pm 1}\,\textup{mod}\,\,p$.
\end{thm}
\begin{rem}
Since $H_2(L_{p,q};\Z)=0$, instead of Poincar\'{e} duality one uses the torsion linking form 
\bgd
H_1(L_{p,q};\Z)\times H_1(L_{p,q};\Z)\to \Q/\Z,
\edd
which classifies these spaces up to homotopy equivalence. Using Reidemeister torsion, defined as 
\bgd
\Delta(L_{p,q}):=(t-1)(t^{q^{-1}}-1)\in \Q[t]/(t^p-1),
\edd
a classification up to PL homeomorphism was given by Reidemeister. Later on, Brody showed that this is a homeomorphism classification. 
\end{rem}
It follows that $L_{7,1}$ and $L_{7,2}$ are homotopy equivalent but not homeomorphic. In fact, even more is true.
\begin{thm}{\bf (Longoni-Salvatore)}\\
The configuration spaces of $L_{7,1}$ and $L_{7,2}$ are not homotopy equivalent.
\end{thm}
\begin{rem}
The proof starts by passing to the 49-sheeted universal covers of the configuration space of two distinct points in $L_{7,j}$. Then these resulting spaces can be distinguished by Massey products. In particular, for $j=1$ the space is formal, being homotopic to $(\vee_6 S^2)\times S^3$ while for $j=2$ there are non-trivial Massey products. This can then be used to distinguish higher configuration spaces.
\end{rem}
\hf\hf We shall denote the $2$-point configuration space of $M$ by $F_2(M)$ as opposed to $C_2(M)$ since in a discussion of chain complexes the latter notation is bound to cause a confusion. Let $M$ be a lens space $L_{p,q}$. We know that the universal cover of $M-\{x_0\}$ is $S^3-\{\Z_p \cdot x_0\}$, which is homotopy equivalent to a bouquet of $(p-1)$ two-spheres. One can imagine replacing the fibre by its universal cover to get a specific $p$-fold cover of $F_2(M)$ that fibres over $M$ with fibre this bouquet. To do this, we use a cross section of the fibration $F_2(M)\to M$ (this always exists since $M$ has a non-vanishing vector field) and replace the fibre $M-\{x_0\}$ by its universal cover which, up to homotopy, is $\vee_{p-1} S^2$. It's necessary to have a cross section for this construction to work since a section gives a consistent choice of base points in the fibre and we replace this with the paths in the fibre based at the base point. Observe that the Hopf fibration $S^3\to S^2$ has no such replacement. Let us denote this new fibration by
\bgd
S^3-\{\Z_p \cdot x_0\}\into F_2(M)'\longrightarrow M.
\edd
We can pull back this fibration using the universal covering map $S^3\to M$ to make a second fibration\ls 
\bgd
S^3-\{\Z_p \cdot x_0\}\into \overline{F_2(M)}\longrightarrow S^3.
\edd
Notice that the total space $\overline{F_2(M)}$ is the universal cover of $F_2(M)$, and the fibre is homotopic to $\vee_{p-1} S^2$. This fibration is determined, up to homotopy, by how the bouquet is mapped to itself along the equator $S^2$ of the base $S^3$, i.e., $\overline{F_2(M)}$ is classified by $\pi_2(\textrm{Aut}_0(\vee_{p-1}S^2))$. Also note that 
\bgd
\overline{F_2(M)}=\{(x,y)\in S^3\times S^3\,|\,x\neq e^{i2\pi k/p}y\}.
\edd
Therefore, for any non-integer $t$ the map $s:x\mapsto (x,e^{i2\pi t/p}x)$ is a bona-fide section. This means that the cohomology ring of $\overline{F_2(M)}$ splits as a tensor product of the cohomology rings of the base and the fibre. 
\begin{eg}{\bf (The lens space $L_{7,1}$)}\\
We think of $S^3$ as the unit quarternions and $\Z_7=\lan \zeta\ran$, the group generated by $\zeta=e^{\frac{2\pi i}{7}}$, acts by left translations. Let $S^3\stackrel{\pi}{\longrightarrow} L_{7,1}$ be the covering map with $\mathbf{1}:=\pi(1)$. For $x,y\in L_{7,1}$ we choose lifts $\overline{x},\overline{y}$ in $S^3$ and define a map 
\bgd
\varphi:F_2(L_{7,1})\longrightarrow L_{7,1}\times \left(L_{7,1}\setminus \{\mathbf{1}\}\right)
\edd
\bgd
\varphi(x,y):=\left(x,\pi((\overline{x})^{-1}\overline{y})\right).
\edd
This is a well defined homeomorphism and thereby provides a homeomorphism 
\bgd
\varphi:\overline{F_2(L_{7,1})}\stackrel{\cong}{\longrightarrow}S^3\times \left(S^3-\Z_7\right). 
\edd
Consequently, $\overline{F_2(L_{7,1})}$ being the product of formal spaces has no Massey products (cf. \ref{MassProd} of \S \ref{Mplay}). \\
\hf\hf If we apply the based loop space functor we get
\bgd
\Omega\varphi:\Omega \overline{F_2(L_{7,1})}\stackrel{\cong}{\longrightarrow}\Omega S^3\times \Omega\left(S^3-\Z_7\right), 
\edd
which is a homeomorphism of $H$-spaces, and it induces an isomorphism of Hopf algebras
\bgd
\Phi:H_\lp^\Q\big(\Omega \overline{F_2(L_{7,1})}\big)\stackrel{\cong}{\longrightarrow} H_\lp^\Q (\Omega S^3)\otimes H_\lp^\Q\big(\Omega (S^3-\Z_7)\big).
\edd
Combining the above with Theorem \ref{MM} (refer \cite{MM65}) we may conclude
\bge
U\Big(\pi_\lp^\Q \big(\Omega \overline{F_2(L_{7,1})}\big)\Big)\stackrel{\cong}{\longrightarrow}  U\big(\pi_\lp^\Q(\Omega S^3)\big)\otimes U\big(\pi_\lp^\Q(\Omega(\vee_6 S^2))\big),
\ede
where $U(\mathfrak{g})$ is the universal enveloping algebra of the Lie algebra $\mathfrak{g}$.\\
\hf\hf To proceed further, we need to know more about the Lie algebra $\pi_\lp^\Q(\Omega(\vee_6 S^2))$. We shall use the following fairly well known result\ls
\begin{lmm}\label{vees2}
The Lie algebra $\pi_\lp^\Q(\Omega(\vee_k S^2))$ is the free Lie algebra generated by $V=\Q^k$ in degree $1$. 
\end{lmm}
{\bf Proof.}\hf We observe that $\vee_k S^2$ is a formal space being the wedge of formal spaces. A minimal model is given by taking $k$ closed generators $x_i$ in degree $2$ and then adding ${k+1 \choose 2}$ generators $\{y_{ij}\}_{1\leq i\leq j\leq k}$ with the relation $dy_{ij}=x_i x_j$. Any higher relations that appear thus far need to be killed appropriately. The terms that need to be killed look like graded commutativity of the Lie bracket dual to the product in this model. In other words, the rank of $\pi_i(\vee_k S^2)$ is given by the dimension of elements of degree $(i-1)$ in the free Lie algebra generated by $V=\Q^k$ of degree $1$. Shifting everything down by $1$, the based loop space of the bouquet with the Samelson product is the free lie algebra on $k$ generators in degree one. $\hfill\square$\\[0.2cm]
Recall Definition \ref{UEA2} of an universal enveloping algebra. It follows that if $\mathfrak{g}$ is a free Lie algebra (on an underlying vector space $V$) then $U\mathfrak{g}$ is simply the free associative algebra on $V$. For the case at hand, $\pi_\lp^\Q(\Omega(\vee_6 S^2))$ is a free Lie algebra, and the universal enveloping algebra is the free associative algebra. In particular, $U\left(\pi_\lp^\Q(\Omega(\vee_6 S^2))\right)$ has no center!
\begin{cor}\label{ueac1}
The center of $U\Big(\pi_\lp^\Q \big(\Omega \overline{F_2(L_{7,1})}\big)\Big)$ is $U(\pi_\lp^\Q(\Omega S^3))\cong \Q[\overline{\alpha}],\,|\overline{\alpha}|=2$. 
\end{cor}
{\bf Proof.}\hf Since $U\left(\pi_\lp^\Q(\Omega(\vee_6 S^2))\right)$ has no center, we conclude that the only central elements arise from $\Omega S^3$. 
Conversely, any element of $U(\pi_\lp^\Q(\Omega S^3))$ is in the centre. Since 
\bgd
\pi_\lp^\Q(\Omega S^3)\cong \Q[\overline{\alpha}]/(\overline{\alpha}^2)
\edd
where $\overline{\alpha}$ generates $\pi_2(\Omega S^3)$, we're done.$\hfill\square$
\end{eg}

%================================ Section 2.2 - Configuration spaces of two points =========================================

\subsection{Configuration spaces of two points II}\label{Conf22}

\hf\hf We deal with the lens space $L_{7,2}$ here.
\begin{eg}\label{L72e}{\bf (The lens space $L_{7,2}$)}\\[0.1cm]
{\it This is a very brief overview of \cite{LS05} and we recall all the relevant Massey products. We're just placing the results of \cite{LS05} in a context applicable for us.}\\[0.1cm]
\hf\hf Let $\zeta=e^{\frac{2\pi i}{7}}$ be a primitive $7$th root of unity. The space $\overline{F_2(L_{7,2})}$ is the complement of the union of the {\it diagonals} in $S^3\times S^3$, i.e.,
\bgd
\overline{F_2(L_{7,2})}=S^3\times S^3-\left(\cup_{k=0}^6 \Delta_k\right),
\edd
where $\Delta_k:=\{(x,\zeta^k x)\,|\,x\in S^3\}$. By the Alexander-Lefschetz duality we have an isomorphism 
\bgd
H^p\big(\overline{F_2(L_{7,2})}\big)\cong H_{6-p}\big(S^3\times S^3,\left(\cup_{k=0}^6 \Delta_k\right)\!\big).
\edd
This identifies the cup product in cohomology with the intersection product in homology. Let $A_k\cong S^3\times [0,1]$ be the submanifold defined by elements of the form
\bgd
((x_1,x_2),(\zeta^{k-1+t}x_1,\zeta^{2(k-1+t)}x_2)),\,\,t\in[0,1],(x_1,x_2)\in S^3.
\edd
\hf\hf Define an action of $\zeta:A_k\to A_{k+1}$ via the map
\bgd
\zeta:((x_1,x_2),(\zeta^{k-1+t}x_1,\zeta^{2(k-1+t)}x_2))\mapsto ((x_1,x_2),(\zeta^{k+t}x_1,\zeta^{2(k+t)}x_2)).
\edd
Since $A_1\trans A_4$ in the interior (as proved in Lemma $4.2$ \cite{LS05}) and this action preserves transversality, we conclude that $A_k\trans A_{k+3}$, where the indices are considered modulo $7$. Define
\bgd
\mathbb{D}_{k,k+3}:=\left\{\big((r,x),(\zeta^{4t+k-1}r,\zeta^{t+2k-2}x)\big)\,|\,(r,x)\in S^3, r^2+|x|^2=1, t\in [0,1]\right\}.
\edd
Notice that $\mathbb{D}_{k,k+3}$ has three boundary components, one each for $r=0, t=0$ and $t=1$. It can be shown that
\bgd
\del_{r=0}\mathbb{D}_{k,k+3}=\left\{\big((0,x),(0,\zeta^{2k-2+t}x)\big)\,|\,|x|=1,t\in[0,1]\right\}=A_k\cap A_{k+3}
\edd
while 
\bgd
\del_{t=0}\mathbb{D}_{k,k+3}\subset \Delta_{k-1},\,\,\,\del_{t=1}\mathbb{D}_{k,k+3}\subset\Delta_{k+3}.
\edd
Similarly, one concludes $\mathbb{D}_{k,k+3}\cap A_{k+5}=\phi$ and
\bgd
A_{k+1}\cap ((1,0)\times S^3)=\mathbb{D}_{k,k+3}\cap A_{k+1}.
\edd
Let us denote the dual (in cohomology) of $A_k$ by $a_k$ and the dual of $(1,0)\times S^3$ by $\alpha$. We can then conclude the following theorem of \cite{LS05}\ls
\begin{thm}{\bf (Longoni-Salvatore)} \label{LSM}\\
The Massey product $\lan a_{k+3},a_k, a_{k+1}+a_{k+5}\ran$ contains the class $a_{k+1}\cup \alpha$. 
\end{thm}
The non-triviality of the Massey product is clear and is essentially due to the addition of $a_{k+5}$ to $a_{k+1}$. Observe that $a_0+\cdots +a_6=0$. This is most transparent when thinking of $\overline{F_2(L_{7,2})}$ as a fibration over $S^3$ with fibre a seven-punctured $S^3$. In other words, the fibre is homotopy equivalent to a wedge of six $2$-spheres and the top cohomology classes of these generate $H^2\big(\overline{F_2(L_{7,2})}\big)$. In fact, with suitable orientations, one can argue that $a_k$ is the fibrewise volume form on the normal sphere bundle of $A_k$ inside $S^3\times S^3$. 
\end{eg}
We shall take up the calculation of the centre of the universal enveloping algebra associated to $\Omega \overline{F_2(L_{7,2})}$ in the next section. The required twisted version will also be analyzed. But before we proceed we need to sort out a few necessary transversality issues.\\
\hf\hf Recall the definition of $\mathbb{D}_{i,i+3}$. We have a similar definition for $\mathbb{D}_{i,j}$. It's birth is necessitated by $A_i\cap A_j$. It can be checked that for $i\neq j$, there are no interior intersection points of $A_i\cap A_j$ unless $j=i+3,i+4$. Therefore, we may set
\bge\label{D}
\mathbb{D}_{i,j}:=\phi\,\,\textup{if}\,\,j\neq i,i+3,i+4.
\ede
One may check the following\ls
\begin{eqnarray}
\label{Dint1}\mathbb{D}_{i,i+3}\cap A_{i+5} & = & \phi\\
\label{Dint2}\mathbb{D}_{i,i+3}\cap A_{i+4} & = & \phi\\
\label{Dint3}\mathbb{D}_{i,i+3}\cap A_{i+1} & = & A_{i+1}\cap (\left(1,0)\times S^3\right)\\
\label{Dint4}\mathbb{D}_{i,i+3}\cap A_{i+2} & = & A_{i+2}\cap (\left(1,0)\times S^3\right).
\end{eqnarray}
The intersections considered on both sides (and as it applies to the equalities as well) are counted modulo boundary. As before, let $a_i$ denote the dual, in cohomology, of $A_i$. One can write all possible Massey products between the $a_i$'s if need be. \\
\hf\hf We switch tracks to discuss a {\it minimal model} for $\overline{F_2(L_{7,2})}$, the universal cover of $F_2(L_{7,2})$. Let $\Lambda_F$ be the minimal model for $\vee_6 S^2$ and $\Lambda(\alpha),d\alpha=0,|\alpha|=3$ be the minimal model for $S^3$. Notice that $\Lambda_F$ starts off with $6$ generators $x_i$ in degree $2$ and then it is formal thereafter. More precisely one needs to {\it kill} $x_i x_j$ by setting $dy_{ij}=x_i x_j$ for $i\leq j$. Similarly, we need $z_{ijk}$ such that
\bgd
dz_{ijk}=x_i y_{jk}-y_{ij}x_k.
\edd
This last equation is very reminiscent of triple Massey products. Indeed, without further terms this says that the Massey products vanish. But our model for $\overline{F_2(L_{7,2})}$ is given by $\Lambda(\alpha)\otimes\Lambda_F$ with the differential $D$ twisted by the Massey product. For example,  we conclude that
\begin{eqnarray}
\label{L72mm1}D(z_{ijk}) & = & x_{i}y_{jk}-y_{ij}x_k\,\,\textup{if}\,\, j\neq i\pm 3\,\textup{and}\, j\neq k\pm  3\\
\label{L72mm2}D(z_{(k+3)k(k+1)}) & = & x_{k+3}y_{k(k+1)}-y_{(k+3)k}x_{k+1}+x_{k+1}\alpha\\
\label{L72mm3}D(z_{(k+3)k(k+5)}) & = & x_{k+3}y_{k(k+5)}-y_{(k+3)k}x_{k+5}.
\end{eqnarray}
One can justify the first equation, for instance, by checking the appropriate intersections, i.e., the homological version of the Massey product and using \eqref{D} and \eqref{Dint1}-\eqref{Dint4}. The other ones are similar. Notice that we're essentially translating the existence of Massey product as the existence of Whitehead products which is tantamount to the non-triviality of the Samelson product on the based loop space. In particular, it follows from the equations above that 
\bgd
D(z_{412}+z_{416})=x_2\alpha +x_4(y_{12}+y_{16})-y_{41}(x_2+x_6).
\edd
The interested reader can check that there are $21$ $y_{ij}$'s and $70$ $z_{ijk}$'s. 

%===================================== Section 2.3 - Universal enveloping algebra ==========================================

\subsection{Pontrjagin product via universal enveloping algebras}

\hf\hf Recall that $\vee_6 S^2\into \overline{F_2(L_{7,i})}\to S^3$ is a fibration of simply connected spaces we may apply the based loop functor to get another fibration
\bgd
\Omega(\vee_6 S^2)\into \Omega\overline{F_2(L_{7,i})}\longrightarrow\Omega S^3.
\edd
Each of the arrows above induce a map of Lie algebras of rational homotopy groups with the Samelson product as the Lie bracket. Further applying the universal enveloping algebra functor we have 
\bge\label{UVW}
\underbrace{U\left(\pi^\Q_\lp(\Omega(\vee_6 S^2))\right)}_{\mathcal{V}}\longrightarrow \underbrace{U\left(\pi^\Q_\lp(\Omega\overline{F_2(L_{7,i})})\right)}_{\mathcal{U}_i}\longrightarrow \underbrace{U\left(\pi^\Q
_\lp(\Omega S^3)\right)}_{\mathcal{W}}.
\ede
By Theorem \ref{MM} of Milnor and Moore \cite{MM65}, the Hurewicz map induces an isomorphism of Hopf algebras
\bgd
\mathcal{H}:U(\pi_\lp^\Q(\Omega X))\stackrel{\cong}{\longrightarrow} \left(H_\lp(\Omega X;\Q),\times\right).
\edd
Using this (or via minimal model description of the original fibration) one checks surjectivity of \eqref{UVW} at $\mathcal{W}$ and injectivity at $\mathcal{V}$. Thus, the objects of interest appear as the middle object in short exact sequences of Hopf algebras 
\bgd
1\to \mathcal{V} \to \mathcal{U}_i \to \mathcal{W}\to 1.
\edd
However, we have already seen that 
\bgd
U\left(\pi^\Q_\lp(\Omega\overline{F_2(L_{7,1})})\right)\cong U\left(\pi^\Q
_\lp(\Omega S^3)\right)\otimes U\left(\pi^\Q_\lp(\Omega(\vee_6 S^2))\right).
\edd
Moreover, it has one central element (up to scaling) in every even degree. On the other hand, the existence of non-trivial Massey products in $F_2(L_{7,2})$ may imply the non-existence of central elements. Our purpose at hand is precisely this. \\
\hf\hf Recall the minimal model for $\overline{F_2(L_{7,2})}$ given by \eqref{L72mm1}-\eqref{L72mm3} from the last section. 
\begin{lmm}
For any simply connected space $X$ let $\Lambda_X$ be a quadratic\footnote{A minimal model for a simply connected space will be called {\it quadratic} if the differential has some non-zero quadratic terms.} minimal model for $X$. Then a minimal model for $\Omega X$ is given by shifting the generators of $\Lambda_X$ down by $1$ and setting the differential to be zero. Moreover, the dual to the Pontrjagin product gives rise to a dg coalgebra structure on $\Lambda_{\Omega X}$. 
\end{lmm}
This is, quite possibly, fairly well known but we sketch a proof for completeness.\\[0.2cm]
{\bf Proof.}\hf The underlying vector space of a minimal model for such a space can be taken to be the rational homotopy groups. Then the quadratic term of the differential in this model is dual to the Whitehead product which is equivalent to the Samelson product on the based loop space. But the Samelson product is the commutator induced by the Pontrjagin product. $\hfill\square$\\[0.2cm]
\hf\hf A minimal model (as a coalgebra) for $\Omega \overline{F_2(L_{7,2})}$ is given by 
\bgd
\Lambda(\overline{\alpha})\otimes \Lambda(\overline{x}_i;\overline{y}_{ij};\overline{z}_{ijk};\cdots), d\equiv 0
\edd
with the coalgebra structure arising from the differential $D$. Here we have used the convention that $\overline{x}$ denotes the element $x$ shifted down in degree by $1$. We shall also interchangeably denote the elements of homotopy groups by $\bar{x}$ etc. as well.
\begin{thm}\label{L72}
There is no element of degree $2$ in the centre of $U\left(\pi^\Q_\lp(\Omega\overline{F_2(L_{7,2})})\right)$.
\end{thm}
{\bf Proof.}\hf It follows from Theorem \ref{LSM} that
\bgd
a_2\cup\alpha \in\lan a_4, a_1, a_2+a_6\ran.
\edd
After appropriately identifying $\overline{x}_i$ with the dual of $a_i$, we conclude that $\overline{\alpha}$ is not central because
\bgd
\overline{x}_2\otimes\overline{\alpha}-\overline{\alpha}\otimes\overline{x}_2=\overline{z}_{412}+\overline{z}_{416}.
\edd
Let $\xi=\lambda\overline{\alpha}+\sum_{i\leq j}\lambda_{ij}\overline{y}_{ij}$ denote an element of degree $2$. If $\xi$ is central then $ [\overline{x}_2,\xi]=0$. However, using the minimal model we have
\begin{eqnarray*}
[\overline{x}_2,\xi] & = & \lambda\overline{z}_{412}+\lambda\overline{z}_{416} + \sum_{i\leq j}\lambda_{ij}(\overline{z}_{2ij}-\overline{z}_{ij2}).
\end{eqnarray*}
The term $\lambda\overline{z}_{416}$ cannot be cancelled by anti-symmetry or by Jacobi identity unless $\lambda=0$. In that case, $\xi$ is a central element in $U\big(\pi_\lp^\Q(\Omega (\vee_6 S^2))\big)$, which being a free associative algebra has no center (see the discussion following Lemma \ref{vees2}). Therefore, $\lambda_{ij}\equiv 0$ and $\xi=0$. $\hfill\square$
\begin{cor}\label{Massey-Pontrjagin}
The Pontrjagin ring $\big(H_\lp(\Omega \overline{F_2(L_{7,2})};\Q),\times\big)$ detects Massey products.
\end{cor}
\hf\hf We would like to prove an analogous statement about the centre of the universal enveloping algebra, but with respect to a {\it twisted} multiplication. It may seem like a stand-alone result with no oversight but we do plan to make crucial use of it in subsequent publication(s) where we show that transversal string topology is {\it not} a homotopy invariant! This has already been emphasized in the introduction. With that being said, we begin with a definition.
\begin{defn}\label{twmul}{\bf (Twisted multiplication)}
Let $\mathcal{A}$ be a (graded) algebra. For $\beta\in\mathcal{A}$ we {\it twist} the multiplication and define
\bgd
x\cdot_{\beta}y:=x\beta y,\,\,\,\,\,x,y\in\mathcal{A}.
\edd
\end{defn}
If $\beta$ has an inverse (and, necessarily, $\mathcal{A}$ has a unit) then we can recover the usual multiplication\ls
\bgd
xy=\beta(\beta^{-1}x\cdot_{\beta}\beta^{-1}y).
\edd
In other words, one {\it left translates} $x$ and $y$ by $\beta^{-1}$ and then uses the twisted multiplication and then {\it right translates} by $\beta$. An element $\xi\in(\mathcal{A},\cdot_\beta)$ is in the {\it centre} if for any $x\in\mathcal{A}$, we have
\bge\label{twcenter}
\xi\beta x=(-1)^{|x||\beta|+|\beta||\xi|+|x||\xi|}x\beta\xi. 
\ede
\begin{defn}
Let $U_{\tau,i}:=U\left(\pi^\Q_\lp(\Omega\overline{F_2(L_{7,i})})\right)$ be the algebra induced by twisting the universal enveloping algebra by the element $\overline{a}_0$. 
\end{defn}
Notice that all powers of $\overline{\alpha}$ are in the centre of $U_{\tau,1}$ since it commutes with everything. On the other hand we have\ls
\begin{thm}\label{L72tw}
There is no element of degree $2$ in the centre of $U_{\tau,2}$. 
\end{thm}
{\bf Proof.}\hf We first show that $\overline{\alpha}$ is not in the centre. We must keep in mind that the notation used in Theorem \ref{LSM} uses $a_i$ and these the same as $x_i$'s used in our model. Moreover, we'll replace $a_0$ by $-a_1-\cdots-a_6$ whenever needed. Let us denote by $\mathcal{Z}$ the following\ls
\begin{eqnarray*}
\mathcal{Z} & := & \overline{\alpha}\otimes \overline{a}_0\otimes 1-1\otimes\overline{a}_0\otimes\overline{\alpha}\\
& = & \sum_{i=1}^6\left(1\otimes\overline{a}_i\otimes\overline{\alpha}-1\otimes\overline{\alpha}\otimes\overline{a}_i\right)\\
& = & \sum_{i=1}^6\left(\overline{a}_i\otimes\overline{\alpha}-\overline{\alpha}\otimes\overline{a}_i\right).
\end{eqnarray*}
Recall the Massey products of Theorem \ref{LSM} as well as those deduced from \eqref{L72mm1}-\eqref{L72mm3}\ls, 
\begin{eqnarray*}
a_2\cup\alpha & \in & \lan a_4, a_1, a_2+a_6\ran\\
a_4\cup\alpha & \in & \lan a_6, a_3, a_4+a_1\ran\\
a_6\cup\alpha & \in & \lan a_1, a_5, a_6+a_3\ran\\
a_3\cup\alpha & \in & \lan a_5, a_2, -a_1-a_2-a_4-a_5-a_6\ran\\
a_5\cup\alpha & \in & \lan -a_1-\cdots -a_6,a_4, a_5+a_2\ran\\
a_1\cup\alpha & \in & \lan a_3, -a_1-\cdots -a_6, a_1+a_5\ran.
\end{eqnarray*}
The above equations, in the language of universal enveloping algebras, should be read as follows. For instance, the fourth equation translates to
\bgd
\overline{z}_{521}+\overline{z}_{522}+\overline{z}_{524}+\overline{z}_{525}+\overline{z}_{526}=\overline{a}_3\otimes\overline{\alpha}-\overline{\alpha}\otimes\overline{a}_3.
\edd
Translating everything relevant to our model we get\ls
\begin{eqnarray}
\label{zee}\mathcal{Z} & = & \overline{z}_{311}+\overline{z}_{315}+\overline{z}_{321}+\overline{z}_{325}+\overline{z}_{331}+\overline{z}_{335}+\overline{z}_{341}+\overline{z}_{345}\\
&   & +\,\overline{z}_{351}+\overline{z}_{355}+\overline{z}_{361}+\overline{z}_{365}-\overline{z}_{412}-\overline{z}_{416}-
\overline{z}_{634}-\overline{z}_{631}\nonumber\\
&   & +\,\overline{z}_{521}+\overline{z}_{522}+\overline{z}_{524}+\overline{z}_{525}+\overline{z}_{526}-\overline{z}_{156}-\overline{z}_{153}\nonumber\\
&   & +\,\overline{z}_{145}+\overline{z}_{142}+\overline{z}_{245}+\overline{z}_{242}+\overline{z}_{345}+\overline{z}_{342}+
\overline{z}_{445}+\overline{z}_{442}\nonumber\\
&   & +\,\overline{z}_{545}+\overline{z}_{542}+\overline{z}_{645}+\overline{z}_{642}\nonumber\\
& \neq & 0.\nonumber
\end{eqnarray}
The easiest way to see why $\mathcal{Z}\neq 0$ is to observe that $\overline{z}_{416}$ doesn't cancel with any other terms. Let
\bgd
\xi=\lambda\overline{\alpha}+\sum_{i\leq j}\lambda_{ij}\overline{y}_{ij}
\edd
denote any central element, with respect to this twisted multiplication, of degree $2$. Then 
\begin{eqnarray*}
\xi\otimes\overline{a}_0\otimes 1-1\otimes\overline{a}_0\otimes\xi & = & \lambda\mathcal{Z}+\sum_{k=1}^6 \sum_{1\leq i\leq j\leq 6}\lambda_{ij}(\overline{z}_{kij}-\overline{z}_{ijk})=0.
\end{eqnarray*}
We shall be making use of the Jacobi identity
\bgd
\overline{z}_{ijk}+\overline{z}_{jki}+\overline{z}_{kij}=0,\,\,\,\,\,i\neq j\neq k.
\edd
We will also use $\overline{z}_{ijk}=-\overline{z}_{kji}$ and $\overline{z}_{iji}=0$ in what follows. 
\begin{eqnarray}
\sum_{k=1}^6 \sum_{1\leq i\leq j\leq 6}\lambda_{ij}(\overline{z}_{ijk}-\overline{z}_{kij}) & = & \sum_{i<j, k\neq i,k\neq j}\lambda_{ij}(\overline{z}_{ijk}-\overline{z}_{kij}) -\sum_{i\neq k}\lambda_{ii}(\overline{z}_{kii}-\overline{z}_{iik})-\sum_{i<j}\lambda_{ij}(\overline{z}_{iij}-\overline{z}_{ijj})\nonumber\\
& = & \sum_{i<j, k\neq i,k\neq j}\lambda_{ij}(\overline{z}_{ijk}-\overline{z}_{kij})+\sum_{i<k}2\lambda_{ii}\overline{z}_{iik}-\sum_{i>k}2\lambda_{ii}\overline{z}_{kii}-\sum_{i<j}\lambda_{ij}(\overline{z}_{iij}-\overline{z}_{ijj})\nonumber\\
\label{cen} & = & \sum_{i<j, k\neq i,k\neq j}\lambda_{ij}(\overline{z}_{ijk}-\overline{z}_{kij})+\sum_{i<j}(2\lambda_{ii}-\lambda_{ij})\overline{z}_{iij}-\sum_{k<l}(2\lambda_{ll}-\lambda_{kl})\overline{z}_{kll}.
\end{eqnarray}
Since $\xi$ is central 
\bgd
\lambda\mathcal{Z}+\sum_{k=1}^6 \sum_{1\leq i\leq j\leq 6}\lambda_{ij}(\overline{z}_{kij}-\overline{z}_{ijk})=0.
\edd
The coefficients of terms of the form $\overline{z}_{ijj}, \overline{z}_{ijj}$ for $i<j$ above should be zero. Using \eqref{zee} and \eqref{cen} we see that
\begin{eqnarray*}
\lambda_{12}=2\lambda_{11} &  & \textup{coefficient of $\overline{z}_{112}$}\\
\lambda_{12}=2\lambda_{22} &  & \textup{coefficient of $\overline{z}_{122}$}\\
\lambda_{14}=2\lambda_{11} &  & \textup{coefficient of $\overline{z}_{114}$}\\
\lambda_{14}=2\lambda_{44} &  & \textup{coefficient of $\overline{z}_{144}$}\\
\lambda_{24}=2\lambda_{22} &  & \textup{coefficient of $\overline{z}_{224}$}\\
\lambda_{24}-\lambda=2\lambda_{44} &  & \textup{coefficient of $\overline{z}_{244}$}.
\end{eqnarray*}
The first four equations imply that $\lambda_{22}=\lambda_{44}$ (they both equal $\lambda_{11}$) which in conjunction with the last two equations imply that $\lambda=0$. This also implies that $\lambda_{ii}=\kappa$ for any $i$ and $\lambda_{ij}=2\kappa$ if $i<j$. Notice that 
\bgd
\xi=\kappa\left(\sum_{i=1}^6 \overline{y}_{ii}+2\sum_{i<j}\overline{y}_{ij}\right)
\edd
now commutes with $\overline{a}_0$. \\
\hf It is enough to show that $\kappa=0$. Consider the twisted commutator\footnote{We are using the Koszul rule of sign where $[\alpha,\beta]_{\tau}:=\alpha\otimes\overline{a}_0\otimes\beta+(-1)^{(|\alpha|+1)(|\beta|+1)}\beta\otimes\overline{a}_0\otimes\alpha$.} of $\xi$ with $\overline{\alpha}$\ls
\begin{eqnarray*}
\xi\otimes\overline{a}_0\otimes\overline{\alpha}-\overline{\alpha}\otimes\overline{a}_0\otimes\xi & = & \overline{a}_0\otimes\xi\otimes\overline{\alpha}-\overline{\alpha}\otimes\overline{a}_0\otimes\xi\\
& = & \overline{a}_0\otimes\overline{\alpha}\otimes\xi-\overline{\alpha}\otimes\overline{a}_0\otimes\xi\\
& = & -\mathcal{Z}\otimes\xi.
\end{eqnarray*}
The result is an element in $U\left(\pi_\lp^\Q(\Omega(\vee_6 S^2))\right)$ which is a (free) tensor algebra. Since $\mathcal{Z}\neq 0$, $\mathcal{Z}\otimes\xi=0$ if and only if $\xi=0$. $\hfill\square$\\[0.2cm]
In the case of $U_{\tau,1}$ the element $\overline{\alpha}$ is clearly in the centre with respect to the twisted multiplication. 
\begin{cor}
The twisted Pontrjagin ring $\big(H_\lp^\Q(\Omega \overline{F_2(L_{7,2})}),\times_{\tau}\big)$ detects Massey products.
\end{cor}

%================================ Section 2.3 - Configuration spaces of n points ===========================================

\subsection{Higher configuration spaces}\label{Confn}

\hf\hf Let $L_{7,j}$ be a lens space obtained by the appropriate action of $\Z_7$ on $S^3$. We consider the universal cover $\overline{F_n(L_{7,j})}$ of $F_n(L_{7,j})$, the configuration space of $n$ points. Our aim is to prove the following\ls
\begin{prpn}\label{lensn}
For $n\geq 2$ the based loop spaces $\Omega \overline{F_n(L_{7,1})}$ and $\Omega \overline{F_n(L_{7,2})}$ are not homotopy equivalent as $H$-spaces.
\end{prpn} 
We proceed to prove this by first showing that central elements exist in every even degree in the case of $L_{7,1}$ while it's absent in those degrees for $L_{7,2}$. Notice that the universal cover of $F_n(L_{7,j})$ is the space of $n$ points in $S^3$ with pairwise disjoint orbit under the $\Z_7$-action. Now consider the map
\bgd
\overline{F_n(L_{7,1})}\longrightarrow S^3\times \mathcal{O}_{n-1},\,\,(x_1,x_2,\cdots,x_n)\mapsto (x_1;x_2x_1^{-1},\cdots,x_n x_1^{-1}),
\edd
where $\mathcal{O}_{n-1}$ is the space of $(n-1)$ points in $S^3-\Z_7$ with pairwise disjoint orbit, or $(n-1)$-{\it orbit configuration space} in short. This is a homeomorphism and induces a map of Hopf algebras
\bge\label{lensn1}
H^\Q_\lp\left(\Omega \overline{F_n(L_{7,1})}\right)\cong H^\Q_\lp(\Omega S^3)\otimes H^\Q_\lp(\Omega \mathcal{O}_{n-1}).
\ede
It is clear that $H^\Q_\lp(\Omega S^3)\cong\Q[\alpha]$, with $|\alpha|=2$, is in the centre. \\
\hf\hf Let us now consider the fibration
\bgd
\mathpzc{F}_n\into \overline{F_n(L_{7,2})}\stackrel{p}{\longrightarrow} \overline{F_{n-1}(L_{7,2})},
\edd
where the map $p$ is given by projecting on the first $n-1$ factors. If we fix $n-1$ points $x_1,x_2,\cdots,x_{n-1}\in S^3$ which have pairwise disjoint orbits then the fibre is
\bgd
\mathpzc{F}_n\cong S^3-\{\Z_7\cdot x_1,\Z_7\cdot x_2,\cdots,\Z_7\cdot x_{n-1}\}\simeq \vee_{7n-8}S^2.
\edd
We now observe that there is a section of this fibration. More precisely, fix a (left-invariant) Riemannian metric on $S^3$ and let $V$ be a (left-invariant) non-vanishing vector field. If $\ep_0$ is the injectivity radius then define
\bge\label{sect}
s(x_1,x_2,\cdots,x_{n-1}):=\left(x_1,\cdots,x_{n-1},\textup{exp}(x_{n-1}+\ep V)\right),
\ede
where $\ep$ is appropriately chosen in terms of $\ep_0$ and the distances between the various orbit points of $x_i$'s. This implies, in particular, that
\bge\label{lensn2}
\pi^\Q_\lp\left(\overline{F_n(L_{7,2})}\right)\cong \pi^\Q_\lp\left(\overline{F_{n-1}(L_{7,2})}\right)\oplus\pi^\Q_\lp (\mathpzc{F}_n).
\ede
Moreover, the left summand above is a Lie subalgebra of the Lie algebra on the left under the Whitehead product. We are now ready to prove our result.\\[0.1cm]
{\bf Proof.}\hf We shall prove our result by induction on $n$ by showing that the Pontrjagin associated to the configuration soace of $n$ points in $L_{7,2}$ has no central element in degree two while it follows from \eqref{lensn1} that the Pontrjagin ring associated to the configuration space of $n$ points in $L_{7,1}$ has central elements in every even degree. \\
\hf We have already dealt with the case $n=2$ (cf. Theorem \ref{L72}). Let us assume we have the statement for $n-1$. Recall the fibration
\bgd
\mathpzc{F}_n\into \overline{F_n(L_{7,2})}\longrightarrow \overline{F_{n-1}(L_{7,2})}
\edd
and the fact that there is a section given by \eqref{sect}. It follows from \eqref{lensn2} that the natural inclusion
\bgd
i:\pi^\Q_\lp (\mathpzc{F}_n)\into \pi^\Q_\lp\left(\overline{F_n(L_{7,2})}\right)
\edd
induces an injective morphism of Hopf algebras
\bgd
i:H^\Q_\lp (\Omega \mathpzc{F}_n) \longrightarrow H^\Q_\lp\left(\Omega \overline{F_n(L_{7,2})}\right).
\edd
Due to the existence of a section, the Pontrjagin ring of $\Omega\mathpzc{F}_n$ is actually an ideal\ls\hspace*{0.05cm}This follows from writing down the minimal model and using \eqref{lensn2}, which implies, for instance, the differential of any element in the fibre cannot have a quadratic term built entirely out of elements from the base. Any degree $2$ element in the Pontrjagin ring associated to $F_n(L_{7,2})$ is given by a degree $3$ element in $\overline{F_n(L_{7,2})}$. Using \eqref{lensn2} let 
\bgd
\alpha\in\pi^\Q_3\left(\overline{F_n(L_{7,2})}\right)=\pi^\Q_3\left(\overline{F_{n-1}(L_{7,2})}\right)\oplus\pi^\Q_3 (\mathpzc{F}_n)
\edd
be written as $\mathpzc{b}\oplus\mathpzc{f}$. If $\mathpzc{b}\neq 0$ then we know by the induction hypothesis that $\mathpzc{b}$ cannot be a central element of degree $2$. Therefore, there exists $\mathpzc{b}'$ such that 
\bgd
\left[\mathpzc{b},\mathpzc{b}'\right]_{\textup{Wh}}\neq 0.
\edd
Moreover, since $\Omega\mathpzc{F}_n$ is an ideal, 
\bgd
\left[\alpha,\mathpzc{b}'\right]_{\textup{Wh}}=\left[\mathpzc{b},\mathpzc{b}'\right]_{\textup{Wh}}\oplus\left[\mathpzc{f},\mathpzc{b}'\right]_{\textup{Wh}}\neq 0.
\edd
Consequently, $\alpha=\mathpzc{f}$ in order to be central. However, the Pontrjagin ring of $\Omega\mathpzc{F}_n$ is the free associative algebra on $7n-8$ generators places in degree one (cf. Lemma \ref{vees2} and the remarks following it). Therefore, $\mathpzc{f}$ can never be central, and there is no central element $\alpha$ as assumed. $\hfill\square$ \\[0.2cm]
\hf\hf We end by observing that there is a short exact sequence of Hopf algebras
\bgd
1\longrightarrow H^\Q_\lp (\Omega \mathpzc{F}_n)\stackrel{i}{\longrightarrow} H^\Q_\lp\left(\Omega \overline{F_n(L_{7,2})}\right)\stackrel{p}{\longrightarrow} H^\Q_\lp\left(\Omega \overline{F_{n-1}(L_{7,2})}\right)\longrightarrow 1.
\edd
Moreover, there is section 
\bgd
s:H^\Q_\lp\left(\Omega \overline{F_{n-1}(L_{7,2})}\right)\longrightarrow H^\Q_\lp\left(\Omega \overline{F_n(L_{7,2})}\right)
\edd
induced by \eqref{sect}. Therefore, both the extreme Hopf algebras in the short exact sequence can be thought of as Hopf subalgebras. We also know that $H^\Q_\lp (\Omega \mathpzc{F}_n)$ is an ideal. Combining these facts it is not hard to see that the {\it twisted\footnote{The twisting is by the same element as in the case for $F_2(L_{7,j})$ since we have an injection $\iota: \pi_2(\Omega \overline{F_2(L_{7,2})})\into \pi_2(\Omega \overline{F_n(L_{7,2})})$.}Pontrjagin rings} associated to $\Omega \overline{F_n(L_{7,2})}$ also do not have central elements in degree $2$. However, it follows from \eqref{lensn1} that there are central elements in every even degree for $\Omega\overline{F_n(L_{7,1})}$. Therefore, the {\it twisted Pontrjagin rings} associated to $\Omega \overline{F_n(L_{7,i})},i=1,2$ are not isomorphic!

\newpage

%============================================ 3 Massey via free loop spaces ================================================

\section{Massey and Pontrjagin products via free loop spaces}\label{MassLM}

\hf\hf Given a minimal model for a simply connected space $X$ (cf. \ref{MinMod} of \S \ref{Mplay}), we can construct a minimal model for the free loop space $LX$. One exploits the existence of a fibration
\bgd
\Omega X\into LX\stackrel{\textup{ev}}{\longrightarrow} X,
\edd
where $\textup{ev}$ is the evaluation of a loop at the base point $1\in S^1$. This model does take into account the twisting present in the fibration. Moreover, one can construct an equivariant model for $LM$ (\cite{CS99}, \S 9) by considering the circle action of rotating the loops. It's possible to embellish these models with further appropriate structures in specific cases to model the loop product, loop bracket and the BV operator on $LM$. However, that is a direction which will lead us into string topology, whose genesis is the paper \cite{CS99} by Chas and Sullivan. This will be pursued elsewhere. \\
\hf\hf In this section we show that the based loop spaces of $S^{2l}\times \Omega S^{2l}$ and $LS^{2l}$, homotopy equivalent as spaces, are not so as $H$-spaces (cf. Theorem \ref{LLSS}). This connection is rather subtle as a similar result fails for $\mathbb{CP}^n$, i.e., the based loop spaces of $\mathbb{CP}^n\times \Omega \mathbb{CP}^n$ and $L\mathbb{CP}^n$ are {\it homotopic} as $H$-spaces (cf. Proposition \ref{OLCP}) if $n\geq 2$. As it turns out, this is strongly correlated with the vanishing of Whitehead products. The story plays out the same in the twisted Pontrjagin ring setup for even spheres versus (complex) projective spaces. In what follows, unless mentioned otherwise, all subsequent models are rational models.

%============================================ 3.1 Minimal models of loop spaces ============================================

\subsection{Models of free loop spaces}\label{stminmod}

\hf\hf We shall begin with the following result (\cite{VPS76}, p. 637) which is also proved in \cite{FHT01} (Example 1, p. 206). 
\begin{thm}\label{SuVP}
If $M$ is a simply connected space with $(\Lambda(x_1,x_2,\cdots),d)$ as its minimal model then the free loop space $LM$ has the minimal model $(\Lambda(x_1,y_1,x_2,y_2,\cdots),\overline{d})$ where $|y_i|=|x_i|-1$. The operator $\overline{d}$ is defined to be $d$ on the $x_i$'s and extended using $\overline{d}s+s\overline{d}=0$, where $s$ is the derivation of $\Lambda(x_1,x_2,\cdots)$ into $\Lambda(x_1,y_1,\cdots)$ defined by $s(x_j)=y_j$.
\end{thm}
Notice that $\Lambda(x_1,x_2,\cdots)$ is a subcomplex of $\Lambda(x_1,y_1,\cdots)$, and the image of $d$ in $\Lambda(x_1,y_1,\cdots)$ is contained in the ideal $I=(x_1,x_2,\cdots)$. Thus, the induced $d$ on $\Lambda(y_j)=\Lambda(x_1,y_1,\cdots)/I$ is zero. This algebraic picture corresponds to the natural fibration
\bgd
\Omega M\into LM\xlongrightarrow{\textrm{ev}} M
\edd
where $(\Lambda(y_j),d\equiv 0)$ can be taken to be a model of $\Omega M$, being a $H$-space. The operator $s$ corresponds to the BV operator on $H_\lp(LM;\Z)$ of rotating loops. \\
\hf\hf Recall that a manifold $M$ is called {\it formal} if there is a quasi-isomorphism between its minimal model $\Lambda_M$ and its cohomology ring, i.e., 
\bgd
\Phi:(\Lambda_M,d)\longrightarrow (H^\lp(M;\Q),0)
\edd
is a quasi-isomorphism. In other words, the rational homotopy theory of the manifold is just a formal consequence of the rational cohomology algebra of $M$. Lie groups and symmetric spaces are known to be formal. Spheres provide a family of examples. It is also known that products and connected sums preserve formality. Moreover, if the Massey product is non-trivial then the manifold is not formal.
\begin{eg}{\bf (K\"ahler manifolds)}\\[0.2cm]
The famous paper \cite{DGMS75} by P. Deligne, P. Griffiths, J. Morgan and D. Sullivan used rational homotopy theory to prove that compact K\"{a}hler manifolds are formal. As a consequence of this rather deep theorem, the complex Grassmanians $G(k,n)$ are formal. In particular, $\proj{C}{n}$ is formal, and this can be verified easily otherwise.  
\end{eg}
\hf\hf We shall work with a specific differential graded algebra $\Lambda(2l)$. This model is just the loop space model of a manifold with a monogenic cohomology ring generated by an element of degree $2l$ and of order $1$. Hence, it covers the case of even dimensional spheres.
\begin{defn}{\bf (The minimal model of $LS^{\textup{even}}$)}\\
For $l\geq 1$ we define $\Lambda(2l)$ to be the differential graded algebra generated by $y_1,x_1,y_2,x_2$ with the differential
\bge
dy_1=0=dx_1,\,\,dx_2=x_1^2,\,\,dy_2=-2x_1y_1.
\ede
The element $x_1$ is of degree $2l$ while $y_1$ is of degree $2l-1$. 
\end{defn}
Before moving on, we emphasize that the algebra above models $LS^{2l}$ in view of Theorem \ref{SuVP} and Example \ref{sph}. Towards calculating the cohomology of the algebra $\Lambda(2l)$, notice that the cochain groups are
\begin{eqnarray*}
C^{\textrm{even}} & = & \Q_\textrm{span}\{y_1x_2 y_2^ax_1^b ,y_2^p x_2^q \}\\
C^{\textrm{odd}} & = & \Q_\textrm{span}\{y_1 y_2^p x_1^q , x_2y_2^ax_1^b \},
\end{eqnarray*}
with the differential given by
\begin{eqnarray}
d(y_1 y_2^p x_1^q ) & = & 0\\
\label{tone}d(y_1x_2 y_2^ax_1^b ) & = & -y_1y_2^a x_1^{b+2}\\
\label{ttne}d(y_2^p x_1^q ) & = & -2py_1y_2^{p-1}x_1^{2+q-1}\\
\label{ne}d(x_2y_2^ax_1^b ) & = & y_2^ax_1^{b+2}-2ay_1x_2y_2^{a-1}x_1^{2+b-1}.
\end{eqnarray}
An element comprising entirely of terms of the form $x_2y_2^ax_1^b$ cannot be closed since the term with the highest $x_1$ exponent survives under the differential. Therefore, the closed odd degree elements are spanned by $y_1 y_2^p x_1^q$ while the image of the even elements only miss $y_1 y_2^p x_1^q$ where $q<1$, i.e., $q=0$. Hence,
\bge\label{lco}
H^{\textrm{odd}}(\Lambda(2l))=\Q_\textrm{span}\{y_1 y_2^p\,|\,p\geq 0\}
\ede
has rank at most $1$. For the even degree elements, 
\bgd
d(y_2^px_1^{q+1}-2py_1x_2y_2^{p-1}x_1^q)=0
\edd
and \eqref{ne} implies that
\bge\label{lce}
H^{\textrm{even}}(\Lambda(2l))=\Q_\textrm{span}\{1,y_2^px_1-2py_1x_2y_2^{p-1}\,|\,p\geq 0\},
\ede
whence the Betti numbers are again $0$ or $1$. In conclusion, we have\ls
\begin{lmm}\label{kl}
The cohomology of the algebra $\Lambda(2l)$ is given by
\begin{eqnarray*}
H^{2j}(\Lambda(2l)) & = & \Q_\textup{span}\{\alpha_{p}\,|\,l+(2l-1)p=j\}\oplus\delta_{j,0}\Q,\\
H^{2j+1}(\Lambda(2l)) & = & \Q_\textup{span}\{y_1 y_2^p \,|\,l+(2l-1)p=j+1\},
\end{eqnarray*}
where
\bgd
\alpha_{p}:=y_2^px_1-2py_1x_2y_2^{p-1}.
\edd
The algebra structure of $H^\lp(\Lambda(2l))$, thought of as the cohomology of $LS^{2l}$ with the wedge product of forms, is trivial.
\end{lmm}
However, the loop space has {\it non-trivial} Massey products, i.e., the loop space is not formal although the base manifold is! This is implied by the following\ls
\begin{thm}\label{MasseyS2l}
The Massey product on the free loop space of even spheres is non-trivial. In particular, with the model $\Lambda(2l)$ and the notations $\alpha_{p}:=y_2^px_1-2py_1x_2y_2^{p-1},\beta_{p}:=y_1y_2^p$, we have 
\begin{eqnarray}
\label{aaa}\left\langle \alpha_{p_1},\alpha_{p_2},\alpha_{p_3}\right\rangle & = & 0\\
\label{bbb}\left\langle \beta_{p_1},\beta_{p_2},\beta_{p_3}\right\rangle & = & 0\\
\label{aab}\left\langle \alpha_{p_1},\alpha_{p_2},\beta_{p_3}\right\rangle & = & \frac{1}{2(p_2+p_3+1)}\alpha_{p+1} \\
\label{bba}\left\langle \beta_{p_1},\beta_{p_2},\alpha_{p_3}\right\rangle & = & \frac{1}{2(p_2+p_3+1)}\beta_{p+1}\\
\label{aba}\left\langle \alpha_{p_1},\beta_{p_2},\alpha_{p_3}\right\rangle & = & \frac{(p_3-p_1)}{2(p_1+p_2+1)(p_2+p_3+1)}\alpha_{p+1}\\
\label{bab}\left\langle \beta_{p_1},\alpha_{p_2},\beta_{p_3}\right\rangle & = & \frac{(p_1-p_3)}{2(p_1+p_2+1)(p_2+p_3+1)}\beta_{p+1},
\end{eqnarray}
where $p=p_1+p_2+p_3$.
\end{thm}
\begin{rem}\label{CPMas}
It is possible to write down a minimal model for the free loop space $L\mathbb{CP}^n$ or $L\mathbb{HP}^n$, where $\mathbb{HP}^n$ is the quarternionic projective space, and prove that Massey products abound in those examples as well! The proof, being essentially identical, is not any more illuminating as the one for even spheres, which we present below.
\end{rem}
{\bf Proof.}\hf Let us denote the sum $p_1+p_2+p_3$ by $p$. Since the cohomology ring is trivial, the triple Massey products are well defined and takes values in $H^\lp(\Lambda(2l))$.\\
\hf \eqref{bbb}\ls\,This is clear.\\
\hf \eqref{aaa}\ls\,When all three are of even degree, then 
\bgd
\alpha_{p_1}\alpha_{p_2}=d(x_2y_2^{p_1+p_2})=dS,\,\,\alpha_{p_2}\alpha_{p_3}=d(x_2y_2^{p_2+p_3})=dT.
\edd
Therefore,
\bgd
\left\langle \alpha_{p_1},\alpha_{p_2},\alpha_{p_3}\right\rangle = \left[S\alpha_{p_3}-\alpha_{p_1}T\right]=0.
\edd
\hf \eqref{aba}\ls\,If $u=\alpha_{p_1},w=\alpha_{p_3},v=\beta_{p_2}$ then 
\begin{eqnarray*}
uv=dS, \,\,S=\frac{-1}{2(p_1+p_2+1)}y_2^{p_1+p_2+1}\\
vw=dT, \,\,T=\frac{-1}{2(p_2+p_3+1)}y_2^{p_2+p_3+1}.
\end{eqnarray*}
A simple manipulation of symbols give us
\bgd
\left\langle u,v,w \right\rangle = \frac{p_3-p_1}{2(p_1+p_2+1)(p_2+p_3+1)}\alpha_{p+1}.
\edd
\hf \eqref{aab}\ls\,Let $u=\alpha_{p_1},v=\alpha_{p_2},w=\beta_{p_3}$ then 
\begin{eqnarray*}
uv=dS,\,\,S=x_2y_2^{p_1+p_2}\\
vw=dT,\,\,T=\frac{-1}{2(p_2+p_3+1)}y_2^{p_2+p_3+1}.
\end{eqnarray*}
The Massey product is given by
\bgd
\left\langle u,v,w\right\rangle = \frac{1}{2(p_2+p_3+1)}\alpha_{p+1}.,
\edd
\hf \eqref{bab}\ls\,Let $u=\beta_{p_1,q_1},v=\alpha_{p_2,q_2}, w=\beta_{p_3,q_3}$. Then 
\begin{eqnarray*}
uv=dS,\,\, S=\frac{-1}{2(p_1+p_2+1)}y_2^{p_1+p_2+1}\\
vw=dT,\,\, T=\frac{-1}{2(p_2+p_3+1)}y_2^{p_2+p_3+1}.
\end{eqnarray*}
The Massey product is given by
\bgd
\left\langle u,v,w\right\rangle = \frac{p_1-p_3}{2(p_1+p_2+1)(p_2+p_3+1)}\beta_{p+1}.
\edd
\hf \eqref{bba}\ls\,Let $u=\beta_{p_1},v=\beta_{p_2}, w=\alpha_{p_3}$. Then $uv=0$ and 
\bgd
vw=dT,\,\,T=\frac{-1}{2(p_2+p_3+1)}y_2^{p_2+p_3+1}.
\edd
The Massey product is given by 
\bgd
\left\langle u,v,w\right\rangle = \frac{1}{2(p_2+p_3+1)}\beta_{p+1},
\edd
This exhausts the list of possibilities for the triple Massey products. $\hfill\square$
\begin{cor}
For the free loop space $LS^{2l}$, the triple Massey products are always defined. In particular, the Massey product for three distinct elements, not of the same parity, is never zero. 
\end{cor}

%============================================ 3.2 Twisted Pontrjagin product ============================================

\subsection{Pontrjagin product on the based loop space of the free loop space}

\hf\hf Recall that we have a fibration of the free loop space over the manifold. In particular, we have
\bgd
\Omega S^{2l}\into LS^{2l}\stackrel{\textup{ev}}{\longrightarrow} S^{2l},
\edd
and there is a section $s:S^{2l}\into LS^{2l}$ given by constant loops. In order to avoid issues of connectedness of $\Omega^2(S^{2l})$, we shall assume $l\geq 2$ unless stated otherwise. We note that Theorem \ref{LLSS} holds for $S^2$ and the proof with the obvious modifications work. We define a map
\bge\label{OSeven}
\Phi: \Omega S^{2l}\times \Omega(\Omega S^{2l})\longrightarrow \Omega (LS^{2l}).
\ede
We belabour on the issue of base points as it is a necessary point! Let $p\in S^{2l}$ and we shall denote by $\mathbbm{p}$ the constant loop based at $p$. We shall also abuse notation and use $\mathbbm{p}$ to denote the base point of $\Omega(\Omega S^{2l})$, now  thought of as the constant loop in $\Omega S^{2l}$ that is based at $\mathbbm{p}$. Given $\gamma\in \Omega  S^{2l},\eta\in \Omega^2(S^{2l})$ we define
\bgd
\Phi(\gamma,\eta):=\eta\ast s(\gamma),
\edd
where $\eta$ and $s(\gamma)$ are to be interpreted as based loops in $LS^{2l}$ and $\ast$ is the concatenation of based loops.\\
\hf\hf It is not hard to check that $\Phi$ is a homotopy equivalence; it induces an isomorphism on homotopy groups and is best seen by identifying the homotopy groups of $\Omega (LS^{2l})$ with a direct sum of (shifted by degree one) homotopy groups of $S^{2l}$ and $\Omega S^{2l}$, and noticing that $\Phi$ is an isomorphism on each summand. Observe that the domain and the range of $\Phi$ in \eqref{OSeven} are formal - based loop spaces being $H$-spaces are formal and formality is preserved under cartesian products. However, $\Phi$ is not an equivalence of $H$-spaces. In fact, there cannot be any $H$-equivalence as the Pontrjagin products on both sides are extremely different! This is a consequence of the following\ls
\begin{thm}\label{LLSS}
For $l\geq 2$ we have the following\ls\\
\textup{(1)} The Pontrjagin rings associated to the spaces $\Omega (LS^{2l})$ and $\Omega S^{2l}\times \Omega^2 S^{2l}$ satisfy\ls\\ 
\hf \textup{(i)} $\big(H_\lp^\Q (\Omega (LS^{2l})),\times\big)$ has no central elements in degree $2l-2$.\\
\hf \textup{(ii)} $\big(H_\lp^\Q (\Omega S^{2l}\times \Omega^2 S^{2l}),\times\big)$ has central elements in degree $2l-2$.\\ 
\textup{(2)} The twisted Pontrjagin rings also exhibit a similar phenomenon\ls\\
\hf \textup{(i)} $\big(H_\lp^\Q (\Omega (LS^{2l})),\times_{\bar{y}_1}\big)$ has no central elements in degree $2l-2$.\\
\hf \textup{(ii)} $\big(H_\lp^\Q (\Omega S^{2l}\times \Omega^2 S^{2l}),\times_{\bar{\beta}_1}\big)$ has central elements in degree $2l-2$.\\
Here $\bar{y}_1$ (or $\bar{\beta}_1$) is a generator of the group $\pi_{2l-2}$, which is $\Z$ for both the spaces above.
\end{thm} 
As will be clear in the course of the proof, a lot more can be said about the exact structure of the Hopf algebras involved but this suffices towards distinguishing the two structures. In order to prove the result, we first develop the minimal models for the two spaces. The result then follows quite easily 
with a little help from algebra. \\
\hf\hf We have seen that a minimal model for $LS^{2l}$ is given by
\bgd
\big(\Lambda(y_1,x_1,y_2,x_2),d\big),\,\,dx_1=0=dy_1,dx_2=x_1^2,dy_2=-2x_1y_1.
\edd
A model of $\Omega(LS^{2l})$ is given by the same algebra as above but by shifting the generators down by one and setting the differential to zero. This follows from the isomorphism
\bgd
\pi_\ast(\Omega(LS^{2l}))\cong\pi_{\ast+1}(LS^{2l})
\edd
and the well known fact in rational homotopy theory that generators in a minimal model (of a simply connected space) can be chosen to be dual to the generators of rational homotopy groups of $LS^{2l}$. We label the dual elements in $\pi^\Q_\ast(LS^{2l})$ as $\bar{x}_1,\bar{x}_2,\bar{y}_1$, and $\bar{y}_2$. Recall that the Whitehead product on the (rational) homotopy groups $\pi^\Q_\lp(LS^{2l})$ of $LS^{2l}$ can be read off from the quadratic part of the differential $d$ in the minimal model. The Whitehead product $[\cdot,\cdot]_{\textup{Wh}}$ on $LS^{2l}$ is precisely the Samelson product on $\pi^\Q_\ast(\Omega(LS^{2l}))$. Therefore, it follows from the model that
\begin{eqnarray*}
\left[\bar{x}_1,\bar{x}_1\right]_{\textup{Wh}} & = & \bar{x}_2\\
\left[\bar{x}_i,\bar{x}_j\right]_{\textup{Wh}} & = & 0\,\,\forall\,\,(i,j)\neq (1,1)\\
\left[\bar{x}_1,\bar{y}_1\right]_{\textup{Wh}} & = & -\frac{\bar{y}_2}{2}\\
\left[\bar{x}_i,\bar{y}_j\right]_{\textup{Wh}} & = & 0\,\,\forall\,\,(i,j)\neq (1,1)\\
\left[\bar{y}_i,\bar{y}_j\right]_{\textup{Wh}} & = & 0\,\,\forall\,\,i,j\in\{1,2\}.\\
\end{eqnarray*}
\hf\hf By Milnor-Moore's theorem \cite{MM65}, the universal enveloping algebra generated by this Lie algebra is the Pontrjagin product on $H_\lp(\Omega(LS^{2l});\Q)$. The degree of the generators get shifted down by one as we are now considering these as elements of $\pi^\Q_\lp(\Omega(LS^{2l}))$. Therefore, the new degree of $\bar{x}_1$ is $2l-1$, that of $\bar{x}_2$ is $4l-2$, that of $\bar{y}_1$ is $2l-2$, and that of $\bar{y}_2$ is $4l-3$. With that being said, the universal enveloping algebra is, by definition, the quotient of the {\it free} tensor algebra $\Lambda$ (generated by $\bar{x}_i$'s and $\bar{y}_i$'s) by the ideal generated by relations of the form
\bgd
a\otimes b-(-1)^{|a||b|}b\otimes a-\left[a,b\right]_{\textup{Wh}}
\edd
for all $a,b\in \Lambda$. If we follow all the degrees and signs posted above, we conclude that the universal enveloping algebra is the free tensor algebra modulo the ideal $I$ generated by the relations 
\begin{eqnarray*}
\bar{x}_1\otimes\bar{x}_1 & = & \frac{1}{2}\bar{x}_2\\
\bar{x}_i\otimes\bar{x}_j & = & (-1)^{|\bar{x}_i||\bar{x}_j|}\bar{x}_j\otimes\bar{x}_i,\,\,(i,j)\neq (1,1)\\
\bar{x}_1\otimes\bar{y}_1 & = & \bar{y}_1\otimes\bar{x}_1-\frac{1}{2}\bar{y}_2\\
\bar{x}_i\otimes\bar{y}_j & = & (-1)^{|\bar{x}_i||\bar{y}_j|}\bar{y}_j\otimes\bar{x}_i,\,\,(i,j)\neq (1,1)\\
\bar{y}_i\otimes\bar{y}_j & = & (-1)^{|\bar{y}_i||\bar{y}_j|}\bar{y}_j\otimes\bar{y}_i,\,\,\forall\,\,i,j\in\{1,2\}.
\end{eqnarray*}
This implies, after a simplication of the algebra, the following\ls
\begin{lmm}\label{LLS}
\textup{(1)} The Pontrjagin ring of $\Omega(LS^{2l})$ is isomorphic to the quotient of the free tensor algebra $T(\bar{x}_1,\bar{x}_2)$, generated by $\bar{x}_1$ and $\bar{y}_1$, by the ideal generated by
\begin{eqnarray*}
\bar{y}_1\otimes\bar{x}_1\otimes\bar{x}_1 & = & \bar{x}_1\otimes\bar{x}_1\otimes \bar{y}_1\\
\bar{x}_1\otimes\bar{y}_1\otimes\bar{y}_1+\bar{y}_1\otimes\bar{y}_1\otimes\bar{x}_1 & = & 2(\bar{y}_1\otimes\bar{x}_1\otimes\bar{y}_1).
\end{eqnarray*}
\textup{(2)} The twisted Pontrjagin ring, twisted by $\bar{y}_1$, of $\Omega(LS^{2l})$ is isomorphic to the underlying vector space as in (1) above. The multiplication is given by $a\times_{\bar{y}_1}b:=a\otimes\bar{y}_1\otimes b$. 
\end{lmm}
{\bf Proof.}\hf It is clear that (2) follows from (1). To prove (1) it suffices to check that all the relations are generated by the relations given above coupled with the requisite commutativity. For instance, $\bar{y}_2\otimes\bar{y}_2=0$ (equivalent to $\bar{y}_2\otimes\bar{y}_2=-\bar{y}_2\otimes\bar{y}_2$) follows from the two relations above and the fact that
\bgd
2(\bar{y}_1\otimes\bar{x}_1 -\bar{x}_1\otimes\bar{y}_1)=\bar{y}_2\\
\edd
The general case is not hard to verify and implies the result.$\hfill\square$\\[0.2cm]
These observations are enough to prove half of the theorem but we will proceed to work on the space $\Omega S^{2l}\times\Omega^2 S^{2l}$.\\
\hf\hf We may decide to work with minimal models as before but there is an easy way out. It is fairly easy to show that there is an isomorphism of Pontrjagin rings
\bge\label{Popro}
H_\lp(\Omega(X\times Y))\cong \big(H_\lp(\Omega X)\otimes H_\lp(\Omega Y)\big).
\ede
This is true with integer coefficients but we shall use $\Q$ for our purposes, and set $X=S^{2l}$ and $Y=\Omega S^{2l}$. We begin with a minimal model for $\Omega S^{2l}$\ls
\bgd
\big(\Lambda(\alpha_1,\alpha_2),d\big),d\equiv 0.
\edd
However, the Lie algabra structure on $\pi^\Q_\lp(\Omega S^{2l})$ (alternatively thought of as $\pi^\Q_{\lp}(S^{2l})$) is given by
\begin{eqnarray*}
\left[\bar{\alpha}_1,\bar{\alpha}_1\right]_{\textup{Wh}} & = & \bar{\alpha}_2\\
\left[\bar{\alpha}_i,\bar{\alpha}_j\right]_{\textup{Wh}} & = & 0\,\,\forall\,\,(i,j)\neq (1,1).
\end{eqnarray*}
Notice that $\bar{\alpha}_1$ has degree $2l-1$ and $\bar{\alpha}_2$ has degree $4l-2$. The universal enveloping algebra associated to this is the free tensor algebra generated by $\bar{\alpha}_1$ as $2\bar{\alpha}_2=\bar{\alpha}_1\otimes\bar{\alpha}_1$. Therefore,
\bge\label{OS}
\big(H_\lp^{\Q}(\Omega S^{2l}),\times\big)\cong T(\bar{\alpha}_1).
\ede
A minimal model for $\Omega^2 S^{2l}$ is given by 
\bgd
\big(\Lambda(\beta_1,\beta_2),d\big)\,\,d\equiv 0,
\edd
where $\beta_1$ has degree $2l-2$ and $\beta_2$ has degree $4l-3$. The Lie algebra structure on $\Omega^2 S^{2l}$ is zero as the differential on $\Omega S^{2l}$ is identically zero. Therefore, the universal enveloping algebra is just the graded commutative algebra on $\bar{\beta}_i$'s, i.e.,
\bge\label{OOS}
\big(H_\lp^{\Q}(\Omega^2 S^{2l}),\times\big)\cong \Lambda(\bar{\beta}_1,\bar{\beta}_2).
\ede
Combining \eqref{OS} and \eqref{OOS} with the observation \eqref{Popro} we conclude that
\bge\label{OSOOS}
\big(H_\lp^{\Q}(\Omega S^{2l}\times \Omega^2 S^{2l}),\times\big)\cong T(\bar{\alpha}_1)\otimes \Lambda(\bar{\beta}_1,\bar{\beta}_2).
\ede
{\bf Proof of Theorem \ref{LLSS}.}\hf Let $l\geq 2$. Part (i) and (ii) of (1) follows easily from Lemma \ref{LLS} and \eqref{OOS}\ls The only possible central element, up to a rational multiple, in $H_{2l-2}(\Omega (LS^{2l}))$ is $\bar{y}_1$ which clearly doesn't commute with $\bar{x}_1$! On the other hand, $\bar{\beta}_1$ commutes with everything due to \eqref{OSOOS}.\\
\hf Part (i) and (ii) of (2) follows since 
\bgd
\bar{y}_1\times_{\bar{y}_1}\bar{x}_1:=\bar{y}_1\otimes \bar{y}_1\otimes \bar{x}_1\neq \bar{x}_1\otimes \bar{y}_1\otimes \bar{y}_1=\bar{x}_1\times_{\bar{y}_1}\bar{y}_1,
\edd
whence $\bar{y}_1$ cannot be central. On the other hand, $\bar{\beta}_1$ is still in the centre as we're twisting by $\bar{\beta}_1$ which is a central element in the untwisted case. $\hfill\square$
\begin{rem}
When $M=S^2$ the based loop space of $LS^2$ is not connected, i.e., $\pi_0(\Omega (LS^2))=\Z$ and the components are indexed by $\pi_2(S^2)$. Similarly, the components of $\Omega^2 S^2$ are also indexed by $\pi_2(S^2)$. We have the usual map
\bgd
\Phi:\Omega S^2 \times \Omega^2 S^2\longrightarrow \Omega (LS^2)
\edd
which is a disjoint union of maps, indexed by $\mathpzc{f}\in\pi_2(S^2)$, given by 
\bgd
\Phi_\mathpzc{f}:\Omega S^2 \times \Omega^2_\mathpzc{f} S^2\longrightarrow \Omega_\mathpzc{f}(LS^2).
\edd
Here $\Omega^2_\mathpzc{f} S^2$ is the connected component of $\mathpzc{f}\in \Omega(\Omega S^2)$, and $\Omega_\mathpzc{f}(LS^2)$ is the connected component of $\mathpzc{f}$ in $\Omega(LS^2)$. Notice that all components of $\Omega^2 S^2$ are homotopy equivalent to $\Omega^2_0S^2$ while all components of $\Omega(LS^2)$ are homotopy equivalent to $\Omega_0(LS^2)$. Therefore, it suffices to distinguish between $\Omega S^2 \times \Omega^2_0 S^2$ and $\Omega_0(LS^2)$ as $H$-spaces. It can be checked that the previous models suffice with the caveat that we take into account the degree zero generator(s).
\end{rem}
\hf\hf We end with a contrasting result. We should remember the underlying meta-principle we're exploiting here\ls\hspace*{0.05cm}If the differential in the minimal model of a space has non-zero quadratic terms then it translates to non-trivial Whitehead products which propagate to Samelson products on the based loop space, and is correlated with Massey products on the space itself. For instance, if we consider the free loop space of $\mathbb{CP}^n$ then we have a minimal model as follows\ls
\bgd
\Lambda(L\mathbb{CP}^n):=(\Lambda(y_1,x_1,y_2,x_2),d),\,\,dy_1=0=dx_1,\,\,dx_2=x_1^{n+1},\,\,dy_2=-(n+1)x_1^ny_1.
\edd
The element $x_1$ is of degree $2$ while $y_1$ is of degree $1$. The Whitehead products vanish if $n\geq 2$ as there are no quadratic terms. However, non-trivial Massey products exist (cf. Remark \ref{CPMas}). The based loop space of either $L\mathbb{CP}^n$ or of $\Omega(\Omega\mathbb{CP}^n)$ is not connected and has $\pi_1(\Omega\mathbb{CP}^n)=\Z$ worth of components. Let us label these components as $\Omega^2_j \mathbb{CP}^n$ and $\Omega_j(L\mathbb{CP}^n)$. Moreover, the fibration 
\bgd
\Omega^2 \mathbb{CP}^n\into \Omega(L\mathbb{CP}^n)\longrightarrow \Omega \mathbb{CP}^n
\edd
is a union of fibrations over $\Omega\mathbb{CP}^n$ where the total space and its corresponding are $\Omega^2_j \mathbb{CP}^n$ and $\Omega_j(L\mathbb{CP}^n)$. The absence of Whitehead products for $L\mathbb{CP}^n$ and $\mathbb{CP}^n$ imply the following\ls
\bge\label{LCP}
\left(H_\lp^{\Q}(\Omega_0(L\mathbb{CP}^n)),\times\right)\cong \Lambda(b,c,d),
\ede
$b,c,d$ are generators of $\pi_i^\Q(L\mathbb{CP}^n)$ for $i=2,2n,2n+1$. Moreover, since $\Omega_j(L\mathbb{CP}^n)$ is homotopic to $\Omega_0(L\mathbb{CP}^n)$ we know $H_\lp^{\Q}(\Omega_j(L\mathbb{CP}^n))$ as an abelian group for any $j$. Finally, the complete ring structure is essentially indexed by $\Z$, i.e., 
\bgd
\times:H_\lp^{\Q}(\Omega_j(L\mathbb{CP}^n))\otimes H_\lp^{\Q}(\Omega_k(L\mathbb{CP}^n))\longrightarrow H_\lp^{\Q}(\Omega_{j+k}(L\mathbb{CP}^n)),
\edd
where the multiplication is the usual multiplication. Written more concisely,
\bge\label{LCP2}
\left(H_\lp^{\Q}(\Omega (L\mathbb{CP}^n)),\times\right)\cong \Lambda(a^{\pm 1},b,c,d),
\ede
where $a,b,c,d$ are generators of $\pi_i^\Q(L\mathbb{CP}^n)$ for $i=1,2,2n,2n+1$, and placed in degrees $0,1,2n-1,2n$ respectively. On the other hand
\begin{eqnarray}
\left(H_\lp^{\Q}(\Omega\mathbb{CP}^n),\times\right) & \cong & U(\pi_{\lp-1}^\Q (\mathbb{CP}^n))\cong\Lambda(\mathpzc{b},\mathpzc{d})\label{OCP}\\
\left(H_\lp^{\Q}(\Omega^2(\mathbb{CP}^n)),\times\right) & \cong & \Lambda(\mathpzc{a}^{\pm 1},\mathpzc{c}) \label{OOCP},
\end{eqnarray}
where $\mathpzc{a},\mathpzc{b},\mathpzc{c},\mathpzc{d}$ have degrees $0,1,2n-1,2n$ respectively. 
\begin{prpn}\label{OLCP}
The natural map $\Phi$, induced by the constant section $s:\mathbb{CP}^n\to L\mathbb{CP}^n$, 
\bgd
\Phi:\Omega \mathbb{CP}^n\times\Omega^2\mathbb{CP}^n\longrightarrow\Omega (L\mathbb{CP}^n)
\edd
is a rational homotopy equivalence of Hopf algebras.
\end{prpn}
{\bf Proof.}\hf The map $\Phi$, as defined above, preserves the Hopf algebra structure as is evidently clear from the obvious isomorphism between \eqref{LCP2} and the tensor product of \eqref{OCP} and \eqref{OOCP}.$\hfill\square$

\bibliographystyle{siam}
\bibliography{Myref_bib.bib}

\vspace*{0.2cm}
\hf {\small D}{\scriptsize EPARTMENT OF }{\small M}{\scriptsize ATHEMATICAL }{\small S}{\scriptsize CIENCES, }{\small B}{\scriptsize INGHAMTON }{\small U}{\scriptsize NIVERSITY, }{\small B}{\scriptsize INGHAMTON, }{\small NY} {\footnotesize 13902-6000}\\
\hf{\it E-mail address} : \texttt{somnath@math.binghamton.edu}

\end{document}

%% file: rgb.tex
\definecolor{AliceBlue}{rgb}{0.94,0.97,1.00}
\definecolor{AntiqueWhite1}{rgb}{1.00,0.94,0.86}
\definecolor{AntiqueWhite2}{rgb}{0.93,0.87,0.80}
\definecolor{AntiqueWhite3}{rgb}{0.80,0.75,0.69}
\definecolor{AntiqueWhite4}{rgb}{0.55,0.51,0.47}
\definecolor{AntiqueWhite}{rgb}{0.98,0.92,0.84}
\definecolor{BlanchedAlmond}{rgb}{1.00,0.92,0.80}
\definecolor{BlueViolet}{rgb}{0.54,0.17,0.89}
\definecolor{CadetBlue1}{rgb}{0.60,0.96,1.00}
\definecolor{CadetBlue2}{rgb}{0.56,0.90,0.93}
\definecolor{CadetBlue3}{rgb}{0.48,0.77,0.80}
\definecolor{CadetBlue4}{rgb}{0.33,0.53,0.55}
\definecolor{CadetBlue}{rgb}{0.37,0.62,0.63}
\definecolor{CornflowerBlue}{rgb}{0.39,0.58,0.93}
\definecolor{DarkBlue}{rgb}{0.00,0.00,0.55}
\definecolor{DarkCyan}{rgb}{0.00,0.55,0.55}
\definecolor{DarkGoldenrod1}{rgb}{1.00,0.73,0.06}
\definecolor{DarkGoldenrod2}{rgb}{0.93,0.68,0.05}
\definecolor{DarkGoldenrod3}{rgb}{0.80,0.58,0.05}
\definecolor{DarkGoldenrod4}{rgb}{0.55,0.40,0.03}
\definecolor{DarkGoldenrod}{rgb}{0.72,0.53,0.04}
\definecolor{DarkGray}{rgb}{0.66,0.66,0.66}
\definecolor{DarkGreen}{rgb}{0.00,0.39,0.00}
\definecolor{DarkGrey}{rgb}{0.66,0.66,0.66}
\definecolor{DarkKhaki}{rgb}{0.74,0.72,0.42}
\definecolor{DarkMagenta}{rgb}{0.55,0.00,0.55}
\definecolor{DarkOliveGreen1}{rgb}{0.79,1.00,0.44}
\definecolor{DarkOliveGreen2}{rgb}{0.74,0.93,0.41}
\definecolor{DarkOliveGreen3}{rgb}{0.64,0.80,0.35}
\definecolor{DarkOliveGreen4}{rgb}{0.43,0.55,0.24}
\definecolor{DarkOliveGreen}{rgb}{0.33,0.42,0.18}
\definecolor{DarkOrange1}{rgb}{1.00,0.50,0.00}
\definecolor{DarkOrange2}{rgb}{0.93,0.46,0.00}
\definecolor{DarkOrange3}{rgb}{0.80,0.40,0.00}
\definecolor{DarkOrange4}{rgb}{0.55,0.27,0.00}
\definecolor{DarkOrange}{rgb}{1.00,0.55,0.00}
\definecolor{DarkOrchid1}{rgb}{0.75,0.24,1.00}
\definecolor{DarkOrchid2}{rgb}{0.70,0.23,0.93}
\definecolor{DarkOrchid3}{rgb}{0.60,0.20,0.80}
\definecolor{DarkOrchid4}{rgb}{0.41,0.13,0.55}
\definecolor{DarkOrchid}{rgb}{0.60,0.20,0.80}
\definecolor{DarkRed}{rgb}{0.55,0.00,0.00}
\definecolor{DarkSalmon}{rgb}{0.91,0.59,0.48}
\definecolor{DarkSeaGreen1}{rgb}{0.76,1.00,0.76}
\definecolor{DarkSeaGreen2}{rgb}{0.71,0.93,0.71}
\definecolor{DarkSeaGreen3}{rgb}{0.61,0.80,0.61}
\definecolor{DarkSeaGreen4}{rgb}{0.41,0.55,0.41}
\definecolor{DarkSeaGreen}{rgb}{0.56,0.74,0.56}
\definecolor{DarkSlateBlue}{rgb}{0.28,0.24,0.55}
\definecolor{DarkSlateGray1}{rgb}{0.59,1.00,1.00}
\definecolor{DarkSlateGray2}{rgb}{0.55,0.93,0.93}
\definecolor{DarkSlateGray3}{rgb}{0.47,0.80,0.80}
\definecolor{DarkSlateGray4}{rgb}{0.32,0.55,0.55}
\definecolor{DarkSlateGray}{rgb}{0.18,0.31,0.31}
\definecolor{DarkSlateGrey}{rgb}{0.18,0.31,0.31}
\definecolor{DarkTurquoise}{rgb}{0.00,0.81,0.82}
\definecolor{DarkViolet}{rgb}{0.58,0.00,0.83}
\definecolor{DeepPink1}{rgb}{1.00,0.08,0.58}
\definecolor{DeepPink2}{rgb}{0.93,0.07,0.54}
\definecolor{DeepPink3}{rgb}{0.80,0.06,0.46}
\definecolor{DeepPink4}{rgb}{0.55,0.04,0.31}
\definecolor{DeepPink}{rgb}{1.00,0.08,0.58}
\definecolor{DeepSkyBlue1}{rgb}{0.00,0.75,1.00}
\definecolor{DeepSkyBlue2}{rgb}{0.00,0.70,0.93}
\definecolor{DeepSkyBlue3}{rgb}{0.00,0.60,0.80}
\definecolor{DeepSkyBlue4}{rgb}{0.00,0.41,0.55}
\definecolor{DeepSkyBlue}{rgb}{0.00,0.75,1.00}
\definecolor{DimGray}{rgb}{0.41,0.41,0.41}
\definecolor{DimGrey}{rgb}{0.41,0.41,0.41}
\definecolor{DodgerBlue1}{rgb}{0.12,0.56,1.00}
\definecolor{DodgerBlue2}{rgb}{0.11,0.53,0.93}
\definecolor{DodgerBlue3}{rgb}{0.09,0.45,0.80}
\definecolor{DodgerBlue4}{rgb}{0.06,0.31,0.55}
\definecolor{DodgerBlue}{rgb}{0.12,0.56,1.00}
\definecolor{FloralWhite}{rgb}{1.00,0.98,0.94}
\definecolor{ForestGreen}{rgb}{0.13,0.55,0.13}
\definecolor{GhostWhite}{rgb}{0.97,0.97,1.00}
\definecolor{GreenYellow}{rgb}{0.68,1.00,0.18}
\definecolor{HotPink1}{rgb}{1.00,0.43,0.71}
\definecolor{HotPink2}{rgb}{0.93,0.42,0.65}
\definecolor{HotPink3}{rgb}{0.80,0.38,0.56}
\definecolor{HotPink4}{rgb}{0.55,0.23,0.38}
\definecolor{HotPink}{rgb}{1.00,0.41,0.71}
\definecolor{IndianRed1}{rgb}{1.00,0.42,0.42}
\definecolor{IndianRed2}{rgb}{0.93,0.39,0.39}
\definecolor{IndianRed3}{rgb}{0.80,0.33,0.33}
\definecolor{IndianRed4}{rgb}{0.55,0.23,0.23}
\definecolor{IndianRed}{rgb}{0.80,0.36,0.36}
\definecolor{LavenderBlush1}{rgb}{1.00,0.94,0.96}
\definecolor{LavenderBlush2}{rgb}{0.93,0.88,0.90}
\definecolor{LavenderBlush3}{rgb}{0.80,0.76,0.77}
\definecolor{LavenderBlush4}{rgb}{0.55,0.51,0.53}
\definecolor{LavenderBlush}{rgb}{1.00,0.94,0.96}
\definecolor{LawnGreen}{rgb}{0.49,0.99,0.00}
\definecolor{LemonChiffon1}{rgb}{1.00,0.98,0.80}
\definecolor{LemonChiffon2}{rgb}{0.93,0.91,0.75}
\definecolor{LemonChiffon3}{rgb}{0.80,0.79,0.65}
\definecolor{LemonChiffon4}{rgb}{0.55,0.54,0.44}
\definecolor{LemonChiffon}{rgb}{1.00,0.98,0.80}
\definecolor{LightBlue1}{rgb}{0.75,0.94,1.00}
\definecolor{LightBlue2}{rgb}{0.70,0.87,0.93}
\definecolor{LightBlue3}{rgb}{0.60,0.75,0.80}
\definecolor{LightBlue4}{rgb}{0.41,0.51,0.55}
\definecolor{LightBlue}{rgb}{0.68,0.85,0.90}
\definecolor{LightCoral}{rgb}{0.94,0.50,0.50}
\definecolor{LightCyan1}{rgb}{0.88,1.00,1.00}
\definecolor{LightCyan2}{rgb}{0.82,0.93,0.93}
\definecolor{LightCyan3}{rgb}{0.71,0.80,0.80}
\definecolor{LightCyan4}{rgb}{0.48,0.55,0.55}
\definecolor{LightCyan}{rgb}{0.88,1.00,1.00}
\definecolor{LightGoldenrod1}{rgb}{1.00,0.93,0.55}
\definecolor{LightGoldenrod2}{rgb}{0.93,0.86,0.51}
\definecolor{LightGoldenrod3}{rgb}{0.80,0.75,0.44}
\definecolor{LightGoldenrod4}{rgb}{0.55,0.51,0.30}
\definecolor{LightGoldenrodYellow}{rgb}{0.98,0.98,0.82}
\definecolor{LightGoldenrod}{rgb}{0.93,0.87,0.51}
\definecolor{LightGray}{rgb}{0.83,0.83,0.83}
\definecolor{LightGreen}{rgb}{0.56,0.93,0.56}
\definecolor{LightGrey}{rgb}{0.83,0.83,0.83}
\definecolor{LightPink1}{rgb}{1.00,0.68,0.73}
\definecolor{LightPink2}{rgb}{0.93,0.64,0.68}
\definecolor{LightPink3}{rgb}{0.80,0.55,0.58}
\definecolor{LightPink4}{rgb}{0.55,0.37,0.40}
\definecolor{LightPink}{rgb}{1.00,0.71,0.76}
\definecolor{LightSalmon1}{rgb}{1.00,0.63,0.48}
\definecolor{LightSalmon2}{rgb}{0.93,0.58,0.45}
\definecolor{LightSalmon3}{rgb}{0.80,0.51,0.38}
\definecolor{LightSalmon4}{rgb}{0.55,0.34,0.26}
\definecolor{LightSalmon}{rgb}{1.00,0.63,0.48}
\definecolor{LightSeaGreen}{rgb}{0.13,0.70,0.67}
\definecolor{LightSkyBlue1}{rgb}{0.69,0.89,1.00}
\definecolor{LightSkyBlue2}{rgb}{0.64,0.83,0.93}
\definecolor{LightSkyBlue3}{rgb}{0.55,0.71,0.80}
\definecolor{LightSkyBlue4}{rgb}{0.38,0.48,0.55}
\definecolor{LightSkyBlue}{rgb}{0.53,0.81,0.98}
\definecolor{LightSlateBlue}{rgb}{0.52,0.44,1.00}
\definecolor{LightSlateGray}{rgb}{0.47,0.53,0.60}
\definecolor{LightSlateGrey}{rgb}{0.47,0.53,0.60}
\definecolor{LightSteelBlue1}{rgb}{0.79,0.88,1.00}
\definecolor{LightSteelBlue2}{rgb}{0.74,0.82,0.93}
\definecolor{LightSteelBlue3}{rgb}{0.64,0.71,0.80}
\definecolor{LightSteelBlue4}{rgb}{0.43,0.48,0.55}
\definecolor{LightSteelBlue}{rgb}{0.69,0.77,0.87}
\definecolor{LightYellow1}{rgb}{1.00,1.00,0.88}
\definecolor{LightYellow2}{rgb}{0.93,0.93,0.82}
\definecolor{LightYellow3}{rgb}{0.80,0.80,0.71}
\definecolor{LightYellow4}{rgb}{0.55,0.55,0.48}
\definecolor{LightYellow}{rgb}{1.00,1.00,0.88}
\definecolor{LimeGreen}{rgb}{0.20,0.80,0.20}
\definecolor{MediumAquamarine}{rgb}{0.40,0.80,0.67}
\definecolor{MediumBlue}{rgb}{0.00,0.00,0.80}
\definecolor{MediumOrchid1}{rgb}{0.88,0.40,1.00}
\definecolor{MediumOrchid2}{rgb}{0.82,0.37,0.93}
\definecolor{MediumOrchid3}{rgb}{0.71,0.32,0.80}
\definecolor{MediumOrchid4}{rgb}{0.48,0.22,0.55}
\definecolor{MediumOrchid}{rgb}{0.73,0.33,0.83}
\definecolor{MediumPurple1}{rgb}{0.67,0.51,1.00}
\definecolor{MediumPurple2}{rgb}{0.62,0.47,0.93}
\definecolor{MediumPurple3}{rgb}{0.54,0.41,0.80}
\definecolor{MediumPurple4}{rgb}{0.36,0.28,0.55}
\definecolor{MediumPurple}{rgb}{0.58,0.44,0.86}
\definecolor{MediumSeaGreen}{rgb}{0.24,0.70,0.44}
\definecolor{MediumSlateBlue}{rgb}{0.48,0.41,0.93}
\definecolor{MediumSpringGreen}{rgb}{0.00,0.98,0.60}
\definecolor{MediumTurquoise}{rgb}{0.28,0.82,0.80}
\definecolor{MediumVioletRed}{rgb}{0.78,0.08,0.52}
\definecolor{MidnightBlue}{rgb}{0.10,0.10,0.44}
\definecolor{MintCream}{rgb}{0.96,1.00,0.98}
\definecolor{MistyRose1}{rgb}{1.00,0.89,0.88}
\definecolor{MistyRose2}{rgb}{0.93,0.84,0.82}
\definecolor{MistyRose3}{rgb}{0.80,0.72,0.71}
\definecolor{MistyRose4}{rgb}{0.55,0.49,0.48}
\definecolor{MistyRose}{rgb}{1.00,0.89,0.88}
\definecolor{NavajoWhite1}{rgb}{1.00,0.87,0.68}
\definecolor{NavajoWhite2}{rgb}{0.93,0.81,0.63}
\definecolor{NavajoWhite3}{rgb}{0.80,0.70,0.55}
\definecolor{NavajoWhite4}{rgb}{0.55,0.47,0.37}
\definecolor{NavajoWhite}{rgb}{1.00,0.87,0.68}
\definecolor{NavyBlue}{rgb}{0.00,0.00,0.50}
\definecolor{OldLace}{rgb}{0.99,0.96,0.90}
\definecolor{OliveDrab1}{rgb}{0.75,1.00,0.24}
\definecolor{OliveDrab2}{rgb}{0.70,0.93,0.23}
\definecolor{OliveDrab3}{rgb}{0.60,0.80,0.20}
\definecolor{OliveDrab4}{rgb}{0.41,0.55,0.13}
\definecolor{OliveDrab}{rgb}{0.42,0.56,0.14}
\definecolor{OrangeRed1}{rgb}{1.00,0.27,0.00}
\definecolor{OrangeRed2}{rgb}{0.93,0.25,0.00}
\definecolor{OrangeRed3}{rgb}{0.80,0.22,0.00}
\definecolor{OrangeRed4}{rgb}{0.55,0.15,0.00}
\definecolor{OrangeRed}{rgb}{1.00,0.27,0.00}
\definecolor{PaleGoldenrod}{rgb}{0.93,0.91,0.67}
\definecolor{PaleGreen1}{rgb}{0.60,1.00,0.60}
\definecolor{PaleGreen2}{rgb}{0.56,0.93,0.56}
\definecolor{PaleGreen3}{rgb}{0.49,0.80,0.49}
\definecolor{PaleGreen4}{rgb}{0.33,0.55,0.33}
\definecolor{PaleGreen}{rgb}{0.60,0.98,0.60}
\definecolor{PaleTurquoise1}{rgb}{0.73,1.00,1.00}
\definecolor{PaleTurquoise2}{rgb}{0.68,0.93,0.93}
\definecolor{PaleTurquoise3}{rgb}{0.59,0.80,0.80}
\definecolor{PaleTurquoise4}{rgb}{0.40,0.55,0.55}
\definecolor{PaleTurquoise}{rgb}{0.69,0.93,0.93}
\definecolor{PaleVioletRed1}{rgb}{1.00,0.51,0.67}
\definecolor{PaleVioletRed2}{rgb}{0.93,0.47,0.62}
\definecolor{PaleVioletRed3}{rgb}{0.80,0.41,0.54}
\definecolor{PaleVioletRed4}{rgb}{0.55,0.28,0.36}
\definecolor{PaleVioletRed}{rgb}{0.86,0.44,0.58}
\definecolor{PapayaWhip}{rgb}{1.00,0.94,0.84}
\definecolor{PeachPuff1}{rgb}{1.00,0.85,0.73}
\definecolor{PeachPuff2}{rgb}{0.93,0.80,0.68}
\definecolor{PeachPuff3}{rgb}{0.80,0.69,0.58}
\definecolor{PeachPuff4}{rgb}{0.55,0.47,0.40}
\definecolor{PeachPuff}{rgb}{1.00,0.85,0.73}
\definecolor{PowderBlue}{rgb}{0.69,0.88,0.90}
\definecolor{RosyBrown1}{rgb}{1.00,0.76,0.76}
\definecolor{RosyBrown2}{rgb}{0.93,0.71,0.71}
\definecolor{RosyBrown3}{rgb}{0.80,0.61,0.61}
\definecolor{RosyBrown4}{rgb}{0.55,0.41,0.41}
\definecolor{RosyBrown}{rgb}{0.74,0.56,0.56}
\definecolor{RoyalBlue1}{rgb}{0.28,0.46,1.00}
\definecolor{RoyalBlue2}{rgb}{0.26,0.43,0.93}
\definecolor{RoyalBlue3}{rgb}{0.23,0.37,0.80}
\definecolor{RoyalBlue4}{rgb}{0.15,0.25,0.55}
\definecolor{RoyalBlue}{rgb}{0.25,0.41,0.88}
\definecolor{SaddleBrown}{rgb}{0.55,0.27,0.07}
\definecolor{SandyBrown}{rgb}{0.96,0.64,0.38}
\definecolor{SeaGreen1}{rgb}{0.33,1.00,0.62}
\definecolor{SeaGreen2}{rgb}{0.31,0.93,0.58}
\definecolor{SeaGreen3}{rgb}{0.26,0.80,0.50}
\definecolor{SeaGreen4}{rgb}{0.18,0.55,0.34}
\definecolor{SeaGreen}{rgb}{0.18,0.55,0.34}
\definecolor{SkyBlue1}{rgb}{0.53,0.81,1.00}
\definecolor{SkyBlue2}{rgb}{0.49,0.75,0.93}
\definecolor{SkyBlue3}{rgb}{0.42,0.65,0.80}
\definecolor{SkyBlue4}{rgb}{0.29,0.44,0.55}
\definecolor{SkyBlue}{rgb}{0.53,0.81,0.92}
\definecolor{SlateBlue1}{rgb}{0.51,0.44,1.00}
\definecolor{SlateBlue2}{rgb}{0.48,0.40,0.93}
\definecolor{SlateBlue3}{rgb}{0.41,0.35,0.80}
\definecolor{SlateBlue4}{rgb}{0.28,0.24,0.55}
\definecolor{SlateBlue}{rgb}{0.42,0.35,0.80}
\definecolor{SlateGray1}{rgb}{0.78,0.89,1.00}
\definecolor{SlateGray2}{rgb}{0.73,0.83,0.93}
\definecolor{SlateGray3}{rgb}{0.62,0.71,0.80}
\definecolor{SlateGray4}{rgb}{0.42,0.48,0.55}
\definecolor{SlateGray}{rgb}{0.44,0.50,0.56}
\definecolor{SlateGrey}{rgb}{0.44,0.50,0.56}
\definecolor{SpringGreen1}{rgb}{0.00,1.00,0.50}
\definecolor{SpringGreen2}{rgb}{0.00,0.93,0.46}
\definecolor{SpringGreen3}{rgb}{0.00,0.80,0.40}
\definecolor{SpringGreen4}{rgb}{0.00,0.55,0.27}
\definecolor{SpringGreen}{rgb}{0.00,1.00,0.50}
\definecolor{SteelBlue1}{rgb}{0.39,0.72,1.00}
\definecolor{SteelBlue2}{rgb}{0.36,0.67,0.93}
\definecolor{SteelBlue3}{rgb}{0.31,0.58,0.80}
\definecolor{SteelBlue4}{rgb}{0.21,0.39,0.55}
\definecolor{SteelBlue}{rgb}{0.27,0.51,0.71}
\definecolor{VioletRed1}{rgb}{1.00,0.24,0.59}
\definecolor{VioletRed2}{rgb}{0.93,0.23,0.55}
\definecolor{VioletRed3}{rgb}{0.80,0.20,0.47}
\definecolor{VioletRed4}{rgb}{0.55,0.13,0.32}
\definecolor{VioletRed}{rgb}{0.82,0.13,0.56}
\definecolor{WhiteSmoke}{rgb}{0.96,0.96,0.96}
\definecolor{YellowGreen}{rgb}{0.60,0.80,0.20}
\definecolor{aliceblue}{rgb}{0.94,0.97,1.00}
\definecolor{antiquewhite}{rgb}{0.98,0.92,0.84}
\definecolor{aquamarine1}{rgb}{0.50,1.00,0.83}
\definecolor{aquamarine2}{rgb}{0.46,0.93,0.78}
\definecolor{aquamarine3}{rgb}{0.40,0.80,0.67}
\definecolor{aquamarine4}{rgb}{0.27,0.55,0.45}
\definecolor{aquamarine}{rgb}{0.50,1.00,0.83}
\definecolor{azure1}{rgb}{0.94,1.00,1.00}
\definecolor{azure2}{rgb}{0.88,0.93,0.93}
\definecolor{azure3}{rgb}{0.76,0.80,0.80}
\definecolor{azure4}{rgb}{0.51,0.55,0.55}
\definecolor{azure}{rgb}{0.94,1.00,1.00}
\definecolor{beige}{rgb}{0.96,0.96,0.86}
\definecolor{bisque1}{rgb}{1.00,0.89,0.77}
\definecolor{bisque2}{rgb}{0.93,0.84,0.72}
\definecolor{bisque3}{rgb}{0.80,0.72,0.62}
\definecolor{bisque4}{rgb}{0.55,0.49,0.42}
\definecolor{bisque}{rgb}{1.00,0.89,0.77}
\definecolor{black}{rgb}{0.00,0.00,0.00}
\definecolor{blanchedalmond}{rgb}{1.00,0.92,0.80}
\definecolor{blue1}{rgb}{0.00,0.00,1.00}
\definecolor{blue2}{rgb}{0.00,0.00,0.93}
\definecolor{blue3}{rgb}{0.00,0.00,0.80}
\definecolor{blue4}{rgb}{0.00,0.00,0.55}
\definecolor{blueviolet}{rgb}{0.54,0.17,0.89}
\definecolor{blue}{rgb}{0.00,0.00,1.00}
\definecolor{brown1}{rgb}{1.00,0.25,0.25}
\definecolor{brown2}{rgb}{0.93,0.23,0.23}
\definecolor{brown3}{rgb}{0.80,0.20,0.20}
\definecolor{brown4}{rgb}{0.55,0.14,0.14}
\definecolor{brown}{rgb}{0.65,0.16,0.16}
\definecolor{burlywood1}{rgb}{1.00,0.83,0.61}
\definecolor{burlywood2}{rgb}{0.93,0.77,0.57}
\definecolor{burlywood3}{rgb}{0.80,0.67,0.49}
\definecolor{burlywood4}{rgb}{0.55,0.45,0.33}
\definecolor{burlywood}{rgb}{0.87,0.72,0.53}
\definecolor{cadetblue}{rgb}{0.37,0.62,0.63}
\definecolor{chartreuse1}{rgb}{0.50,1.00,0.00}
\definecolor{chartreuse2}{rgb}{0.46,0.93,0.00}
\definecolor{chartreuse3}{rgb}{0.40,0.80,0.00}
\definecolor{chartreuse4}{rgb}{0.27,0.55,0.00}
\definecolor{chartreuse}{rgb}{0.50,1.00,0.00}
\definecolor{chocolate1}{rgb}{1.00,0.50,0.14}
\definecolor{chocolate2}{rgb}{0.93,0.46,0.13}
\definecolor{chocolate3}{rgb}{0.80,0.40,0.11}
\definecolor{chocolate4}{rgb}{0.55,0.27,0.07}
\definecolor{chocolate}{rgb}{0.82,0.41,0.12}
\definecolor{coral1}{rgb}{1.00,0.45,0.34}
\definecolor{coral2}{rgb}{0.93,0.42,0.31}
\definecolor{coral3}{rgb}{0.80,0.36,0.27}
\definecolor{coral4}{rgb}{0.55,0.24,0.18}
\definecolor{coral}{rgb}{1.00,0.50,0.31}
\definecolor{cornflowerblue}{rgb}{0.39,0.58,0.93}
\definecolor{cornsilk1}{rgb}{1.00,0.97,0.86}
\definecolor{cornsilk2}{rgb}{0.93,0.91,0.80}
\definecolor{cornsilk3}{rgb}{0.80,0.78,0.69}
\definecolor{cornsilk4}{rgb}{0.55,0.53,0.47}
\definecolor{cornsilk}{rgb}{1.00,0.97,0.86}
\definecolor{cyan1}{rgb}{0.00,1.00,1.00}
\definecolor{cyan2}{rgb}{0.00,0.93,0.93}
\definecolor{cyan3}{rgb}{0.00,0.80,0.80}
\definecolor{cyan4}{rgb}{0.00,0.55,0.55}
\definecolor{cyan}{rgb}{0.00,1.00,1.00}
\definecolor{darkblue}{rgb}{0.00,0.00,0.55}
\definecolor{darkcyan}{rgb}{0.00,0.55,0.55}
\definecolor{darkgoldenrod}{rgb}{0.72,0.53,0.04}
\definecolor{darkgray}{rgb}{0.66,0.66,0.66}
\definecolor{darkgreen}{rgb}{0.00,0.39,0.00}
\definecolor{darkgrey}{rgb}{0.66,0.66,0.66}
\definecolor{darkkhaki}{rgb}{0.74,0.72,0.42}
\definecolor{darkmagenta}{rgb}{0.55,0.00,0.55}
\definecolor{darkolive}{rgb}{0.33,0.42,0.18}
\definecolor{darkorange}{rgb}{1.00,0.55,0.00}
\definecolor{darkorchid}{rgb}{0.60,0.20,0.80}
\definecolor{darkred}{rgb}{0.55,0.00,0.00}
\definecolor{darksalmon}{rgb}{0.91,0.59,0.48}
\definecolor{darksea}{rgb}{0.56,0.74,0.56}
\definecolor{darkslate}{rgb}{0.18,0.31,0.31}
\definecolor{darkslate}{rgb}{0.18,0.31,0.31}
\definecolor{darkslate}{rgb}{0.28,0.24,0.55}
\definecolor{darkturquoise}{rgb}{0.00,0.81,0.82}
\definecolor{darkviolet}{rgb}{0.58,0.00,0.83}
\definecolor{deeppink}{rgb}{1.00,0.08,0.58}
\definecolor{deepsky}{rgb}{0.00,0.75,1.00}
\definecolor{dimgray}{rgb}{0.41,0.41,0.41}
\definecolor{dimgrey}{rgb}{0.41,0.41,0.41}
\definecolor{dodgerblue}{rgb}{0.12,0.56,1.00}
\definecolor{firebrick1}{rgb}{1.00,0.19,0.19}
\definecolor{firebrick2}{rgb}{0.93,0.17,0.17}
\definecolor{firebrick3}{rgb}{0.80,0.15,0.15}
\definecolor{firebrick4}{rgb}{0.55,0.10,0.10}
\definecolor{firebrick}{rgb}{0.70,0.13,0.13}
\definecolor{floralwhite}{rgb}{1.00,0.98,0.94}
\definecolor{forestgreen}{rgb}{0.13,0.55,0.13}
\definecolor{gainsboro}{rgb}{0.86,0.86,0.86}
\definecolor{ghostwhite}{rgb}{0.97,0.97,1.00}
\definecolor{gold1}{rgb}{1.00,0.84,0.00}
\definecolor{gold2}{rgb}{0.93,0.79,0.00}
\definecolor{gold3}{rgb}{0.80,0.68,0.00}
\definecolor{gold4}{rgb}{0.55,0.46,0.00}
\definecolor{goldenrod1}{rgb}{1.00,0.76,0.15}
\definecolor{goldenrod2}{rgb}{0.93,0.71,0.13}
\definecolor{goldenrod3}{rgb}{0.80,0.61,0.11}
\definecolor{goldenrod4}{rgb}{0.55,0.41,0.08}
\definecolor{goldenrod}{rgb}{0.85,0.65,0.13}
\definecolor{gold}{rgb}{1.00,0.84,0.00}
\definecolor{gray0}{rgb}{0.00,0.00,0.00}
\definecolor{gray100}{rgb}{1.00,1.00,1.00}
\definecolor{gray10}{rgb}{0.10,0.10,0.10}
\definecolor{gray11}{rgb}{0.11,0.11,0.11}
\definecolor{gray12}{rgb}{0.12,0.12,0.12}
\definecolor{gray13}{rgb}{0.13,0.13,0.13}
\definecolor{gray14}{rgb}{0.14,0.14,0.14}
\definecolor{gray15}{rgb}{0.15,0.15,0.15}
\definecolor{gray16}{rgb}{0.16,0.16,0.16}
\definecolor{gray17}{rgb}{0.17,0.17,0.17}
\definecolor{gray18}{rgb}{0.18,0.18,0.18}
\definecolor{gray19}{rgb}{0.19,0.19,0.19}
\definecolor{gray1}{rgb}{0.01,0.01,0.01}
\definecolor{gray20}{rgb}{0.20,0.20,0.20}
\definecolor{gray21}{rgb}{0.21,0.21,0.21}
\definecolor{gray22}{rgb}{0.22,0.22,0.22}
\definecolor{gray23}{rgb}{0.23,0.23,0.23}
\definecolor{gray24}{rgb}{0.24,0.24,0.24}
\definecolor{gray25}{rgb}{0.25,0.25,0.25}
\definecolor{gray26}{rgb}{0.26,0.26,0.26}
\definecolor{gray27}{rgb}{0.27,0.27,0.27}
\definecolor{gray28}{rgb}{0.28,0.28,0.28}
\definecolor{gray29}{rgb}{0.29,0.29,0.29}
\definecolor{gray2}{rgb}{0.02,0.02,0.02}
\definecolor{gray30}{rgb}{0.30,0.30,0.30}
\definecolor{gray31}{rgb}{0.31,0.31,0.31}
\definecolor{gray32}{rgb}{0.32,0.32,0.32}
\definecolor{gray33}{rgb}{0.33,0.33,0.33}
\definecolor{gray34}{rgb}{0.34,0.34,0.34}
\definecolor{gray35}{rgb}{0.35,0.35,0.35}
\definecolor{gray36}{rgb}{0.36,0.36,0.36}
\definecolor{gray37}{rgb}{0.37,0.37,0.37}
\definecolor{gray38}{rgb}{0.38,0.38,0.38}
\definecolor{gray39}{rgb}{0.39,0.39,0.39}
\definecolor{gray3}{rgb}{0.03,0.03,0.03}
\definecolor{gray40}{rgb}{0.40,0.40,0.40}
\definecolor{gray41}{rgb}{0.41,0.41,0.41}
\definecolor{gray42}{rgb}{0.42,0.42,0.42}
\definecolor{gray43}{rgb}{0.43,0.43,0.43}
\definecolor{gray44}{rgb}{0.44,0.44,0.44}
\definecolor{gray45}{rgb}{0.45,0.45,0.45}
\definecolor{gray46}{rgb}{0.46,0.46,0.46}
\definecolor{gray47}{rgb}{0.47,0.47,0.47}
\definecolor{gray48}{rgb}{0.48,0.48,0.48}
\definecolor{gray49}{rgb}{0.49,0.49,0.49}
\definecolor{gray4}{rgb}{0.04,0.04,0.04}
\definecolor{gray50}{rgb}{0.50,0.50,0.50}
\definecolor{gray51}{rgb}{0.51,0.51,0.51}
\definecolor{gray52}{rgb}{0.52,0.52,0.52}
\definecolor{gray53}{rgb}{0.53,0.53,0.53}
\definecolor{gray54}{rgb}{0.54,0.54,0.54}
\definecolor{gray55}{rgb}{0.55,0.55,0.55}
\definecolor{gray56}{rgb}{0.56,0.56,0.56}
\definecolor{gray57}{rgb}{0.57,0.57,0.57}
\definecolor{gray58}{rgb}{0.58,0.58,0.58}
\definecolor{gray59}{rgb}{0.59,0.59,0.59}
\definecolor{gray5}{rgb}{0.05,0.05,0.05}
\definecolor{gray60}{rgb}{0.60,0.60,0.60}
\definecolor{gray61}{rgb}{0.61,0.61,0.61}
\definecolor{gray62}{rgb}{0.62,0.62,0.62}
\definecolor{gray63}{rgb}{0.63,0.63,0.63}
\definecolor{gray64}{rgb}{0.64,0.64,0.64}
\definecolor{gray65}{rgb}{0.65,0.65,0.65}
\definecolor{gray66}{rgb}{0.66,0.66,0.66}
\definecolor{gray67}{rgb}{0.67,0.67,0.67}
\definecolor{gray68}{rgb}{0.68,0.68,0.68}
\definecolor{gray69}{rgb}{0.69,0.69,0.69}
\definecolor{gray6}{rgb}{0.06,0.06,0.06}
\definecolor{gray70}{rgb}{0.70,0.70,0.70}
\definecolor{gray71}{rgb}{0.71,0.71,0.71}
\definecolor{gray72}{rgb}{0.72,0.72,0.72}
\definecolor{gray73}{rgb}{0.73,0.73,0.73}
\definecolor{gray74}{rgb}{0.74,0.74,0.74}
\definecolor{gray75}{rgb}{0.75,0.75,0.75}
\definecolor{gray76}{rgb}{0.76,0.76,0.76}
\definecolor{gray77}{rgb}{0.77,0.77,0.77}
\definecolor{gray78}{rgb}{0.78,0.78,0.78}
\definecolor{gray79}{rgb}{0.79,0.79,0.79}
\definecolor{gray7}{rgb}{0.07,0.07,0.07}
\definecolor{gray80}{rgb}{0.80,0.80,0.80}
\definecolor{gray81}{rgb}{0.81,0.81,0.81}
\definecolor{gray82}{rgb}{0.82,0.82,0.82}
\definecolor{gray83}{rgb}{0.83,0.83,0.83}
\definecolor{gray84}{rgb}{0.84,0.84,0.84}
\definecolor{gray85}{rgb}{0.85,0.85,0.85}
\definecolor{gray86}{rgb}{0.86,0.86,0.86}
\definecolor{gray87}{rgb}{0.87,0.87,0.87}
\definecolor{gray88}{rgb}{0.88,0.88,0.88}
\definecolor{gray89}{rgb}{0.89,0.89,0.89}
\definecolor{gray8}{rgb}{0.08,0.08,0.08}
\definecolor{gray90}{rgb}{0.90,0.90,0.90}
\definecolor{gray91}{rgb}{0.91,0.91,0.91}
\definecolor{gray92}{rgb}{0.92,0.92,0.92}
\definecolor{gray93}{rgb}{0.93,0.93,0.93}
\definecolor{gray94}{rgb}{0.94,0.94,0.94}
\definecolor{gray95}{rgb}{0.95,0.95,0.95}
\definecolor{gray96}{rgb}{0.96,0.96,0.96}
\definecolor{gray97}{rgb}{0.97,0.97,0.97}
\definecolor{gray98}{rgb}{0.98,0.98,0.98}
\definecolor{gray99}{rgb}{0.99,0.99,0.99}
\definecolor{gray9}{rgb}{0.09,0.09,0.09}
\definecolor{gray}{rgb}{0.75,0.75,0.75}
\definecolor{green1}{rgb}{0.00,1.00,0.00}
\definecolor{green2}{rgb}{0.00,0.93,0.00}
\definecolor{green3}{rgb}{0.00,0.80,0.00}
\definecolor{green4}{rgb}{0.00,0.55,0.00}
\definecolor{greenyellow}{rgb}{0.68,1.00,0.18}
\definecolor{green}{rgb}{0.00,1.00,0.00}
\definecolor{grey0}{rgb}{0.00,0.00,0.00}
\definecolor{grey100}{rgb}{1.00,1.00,1.00}
\definecolor{grey10}{rgb}{0.10,0.10,0.10}
\definecolor{grey11}{rgb}{0.11,0.11,0.11}
\definecolor{grey12}{rgb}{0.12,0.12,0.12}
\definecolor{grey13}{rgb}{0.13,0.13,0.13}
\definecolor{grey14}{rgb}{0.14,0.14,0.14}
\definecolor{grey15}{rgb}{0.15,0.15,0.15}
\definecolor{grey16}{rgb}{0.16,0.16,0.16}
\definecolor{grey17}{rgb}{0.17,0.17,0.17}
\definecolor{grey18}{rgb}{0.18,0.18,0.18}
\definecolor{grey19}{rgb}{0.19,0.19,0.19}
\definecolor{grey1}{rgb}{0.01,0.01,0.01}
\definecolor{grey20}{rgb}{0.20,0.20,0.20}
\definecolor{grey21}{rgb}{0.21,0.21,0.21}
\definecolor{grey22}{rgb}{0.22,0.22,0.22}
\definecolor{grey23}{rgb}{0.23,0.23,0.23}
\definecolor{grey24}{rgb}{0.24,0.24,0.24}
\definecolor{grey25}{rgb}{0.25,0.25,0.25}
\definecolor{grey26}{rgb}{0.26,0.26,0.26}
\definecolor{grey27}{rgb}{0.27,0.27,0.27}
\definecolor{grey28}{rgb}{0.28,0.28,0.28}
\definecolor{grey29}{rgb}{0.29,0.29,0.29}
\definecolor{grey2}{rgb}{0.02,0.02,0.02}
\definecolor{grey30}{rgb}{0.30,0.30,0.30}
\definecolor{grey31}{rgb}{0.31,0.31,0.31}
\definecolor{grey32}{rgb}{0.32,0.32,0.32}
\definecolor{grey33}{rgb}{0.33,0.33,0.33}
\definecolor{grey34}{rgb}{0.34,0.34,0.34}
\definecolor{grey35}{rgb}{0.35,0.35,0.35}
\definecolor{grey36}{rgb}{0.36,0.36,0.36}
\definecolor{grey37}{rgb}{0.37,0.37,0.37}
\definecolor{grey38}{rgb}{0.38,0.38,0.38}
\definecolor{grey39}{rgb}{0.39,0.39,0.39}
\definecolor{grey3}{rgb}{0.03,0.03,0.03}
\definecolor{grey40}{rgb}{0.40,0.40,0.40}
\definecolor{grey41}{rgb}{0.41,0.41,0.41}
\definecolor{grey42}{rgb}{0.42,0.42,0.42}
\definecolor{grey43}{rgb}{0.43,0.43,0.43}
\definecolor{grey44}{rgb}{0.44,0.44,0.44}
\definecolor{grey45}{rgb}{0.45,0.45,0.45}
\definecolor{grey46}{rgb}{0.46,0.46,0.46}
\definecolor{grey47}{rgb}{0.47,0.47,0.47}
\definecolor{grey48}{rgb}{0.48,0.48,0.48}
\definecolor{grey49}{rgb}{0.49,0.49,0.49}
\definecolor{grey4}{rgb}{0.04,0.04,0.04}
\definecolor{grey50}{rgb}{0.50,0.50,0.50}
\definecolor{grey51}{rgb}{0.51,0.51,0.51}
\definecolor{grey52}{rgb}{0.52,0.52,0.52}
\definecolor{grey53}{rgb}{0.53,0.53,0.53}
\definecolor{grey54}{rgb}{0.54,0.54,0.54}
\definecolor{grey55}{rgb}{0.55,0.55,0.55}
\definecolor{grey56}{rgb}{0.56,0.56,0.56}
\definecolor{grey57}{rgb}{0.57,0.57,0.57}
\definecolor{grey58}{rgb}{0.58,0.58,0.58}
\definecolor{grey59}{rgb}{0.59,0.59,0.59}
\definecolor{grey5}{rgb}{0.05,0.05,0.05}
\definecolor{grey60}{rgb}{0.60,0.60,0.60}
\definecolor{grey61}{rgb}{0.61,0.61,0.61}
\definecolor{grey62}{rgb}{0.62,0.62,0.62}
\definecolor{grey63}{rgb}{0.63,0.63,0.63}
\definecolor{grey64}{rgb}{0.64,0.64,0.64}
\definecolor{grey65}{rgb}{0.65,0.65,0.65}
\definecolor{grey66}{rgb}{0.66,0.66,0.66}
\definecolor{grey67}{rgb}{0.67,0.67,0.67}
\definecolor{grey68}{rgb}{0.68,0.68,0.68}
\definecolor{grey69}{rgb}{0.69,0.69,0.69}
\definecolor{grey6}{rgb}{0.06,0.06,0.06}
\definecolor{grey70}{rgb}{0.70,0.70,0.70}
\definecolor{grey71}{rgb}{0.71,0.71,0.71}
\definecolor{grey72}{rgb}{0.72,0.72,0.72}
\definecolor{grey73}{rgb}{0.73,0.73,0.73}
\definecolor{grey74}{rgb}{0.74,0.74,0.74}
\definecolor{grey75}{rgb}{0.75,0.75,0.75}
\definecolor{grey76}{rgb}{0.76,0.76,0.76}
\definecolor{grey77}{rgb}{0.77,0.77,0.77}
\definecolor{grey78}{rgb}{0.78,0.78,0.78}
\definecolor{grey79}{rgb}{0.79,0.79,0.79}
\definecolor{grey7}{rgb}{0.07,0.07,0.07}
\definecolor{grey80}{rgb}{0.80,0.80,0.80}
\definecolor{grey81}{rgb}{0.81,0.81,0.81}
\definecolor{grey82}{rgb}{0.82,0.82,0.82}
\definecolor{grey83}{rgb}{0.83,0.83,0.83}
\definecolor{grey84}{rgb}{0.84,0.84,0.84}
\definecolor{grey85}{rgb}{0.85,0.85,0.85}
\definecolor{grey86}{rgb}{0.86,0.86,0.86}
\definecolor{grey87}{rgb}{0.87,0.87,0.87}
\definecolor{grey88}{rgb}{0.88,0.88,0.88}
\definecolor{grey89}{rgb}{0.89,0.89,0.89}
\definecolor{grey8}{rgb}{0.08,0.08,0.08}
\definecolor{grey90}{rgb}{0.90,0.90,0.90}
\definecolor{grey91}{rgb}{0.91,0.91,0.91}
\definecolor{grey92}{rgb}{0.92,0.92,0.92}
\definecolor{grey93}{rgb}{0.93,0.93,0.93}
\definecolor{grey94}{rgb}{0.94,0.94,0.94}
\definecolor{grey95}{rgb}{0.95,0.95,0.95}
\definecolor{grey96}{rgb}{0.96,0.96,0.96}
\definecolor{grey97}{rgb}{0.97,0.97,0.97}
\definecolor{grey98}{rgb}{0.98,0.98,0.98}
\definecolor{grey99}{rgb}{0.99,0.99,0.99}
\definecolor{grey9}{rgb}{0.09,0.09,0.09}
\definecolor{grey}{rgb}{0.75,0.75,0.75}
\definecolor{honeydew1}{rgb}{0.94,1.00,0.94}
\definecolor{honeydew2}{rgb}{0.88,0.93,0.88}
\definecolor{honeydew3}{rgb}{0.76,0.80,0.76}
\definecolor{honeydew4}{rgb}{0.51,0.55,0.51}
\definecolor{honeydew}{rgb}{0.94,1.00,0.94}
\definecolor{hotpink}{rgb}{1.00,0.41,0.71}
\definecolor{indianred}{rgb}{0.80,0.36,0.36}
\definecolor{ivory1}{rgb}{1.00,1.00,0.94}
\definecolor{ivory2}{rgb}{0.93,0.93,0.88}
\definecolor{ivory3}{rgb}{0.80,0.80,0.76}
\definecolor{ivory4}{rgb}{0.55,0.55,0.51}
\definecolor{ivory}{rgb}{1.00,1.00,0.94}
\definecolor{khaki1}{rgb}{1.00,0.96,0.56}
\definecolor{khaki2}{rgb}{0.93,0.90,0.52}
\definecolor{khaki3}{rgb}{0.80,0.78,0.45}
\definecolor{khaki4}{rgb}{0.55,0.53,0.31}
\definecolor{khaki}{rgb}{0.94,0.90,0.55}
\definecolor{lavenderblush}{rgb}{1.00,0.94,0.96}
\definecolor{lavender}{rgb}{0.90,0.90,0.98}
\definecolor{lawngreen}{rgb}{0.49,0.99,0.00}
\definecolor{lemonchiffon}{rgb}{1.00,0.98,0.80}
\definecolor{lightblue}{rgb}{0.68,0.85,0.90}
\definecolor{lightcoral}{rgb}{0.94,0.50,0.50}
\definecolor{lightcyan}{rgb}{0.88,1.00,1.00}
\definecolor{lightgoldenrod}{rgb}{0.93,0.87,0.51}
\definecolor{lightgoldenrod}{rgb}{0.98,0.98,0.82}
\definecolor{lightgray}{rgb}{0.83,0.83,0.83}
\definecolor{lightgreen}{rgb}{0.56,0.93,0.56}
\definecolor{lightgrey}{rgb}{0.83,0.83,0.83}
\definecolor{lightpink}{rgb}{1.00,0.71,0.76}
\definecolor{lightsalmon}{rgb}{1.00,0.63,0.48}
\definecolor{lightsea}{rgb}{0.13,0.70,0.67}
\definecolor{lightsky}{rgb}{0.53,0.81,0.98}
\definecolor{lightslate}{rgb}{0.47,0.53,0.60}
\definecolor{lightslate}{rgb}{0.47,0.53,0.60}
\definecolor{lightslate}{rgb}{0.52,0.44,1.00}
\definecolor{lightsteel}{rgb}{0.69,0.77,0.87}
\definecolor{lightyellow}{rgb}{1.00,1.00,0.88}
\definecolor{limegreen}{rgb}{0.20,0.80,0.20}
\definecolor{linen}{rgb}{0.98,0.94,0.90}
\definecolor{magenta1}{rgb}{1.00,0.00,1.00}
\definecolor{magenta2}{rgb}{0.93,0.00,0.93}
\definecolor{magenta3}{rgb}{0.80,0.00,0.80}
\definecolor{magenta4}{rgb}{0.55,0.00,0.55}
\definecolor{magenta}{rgb}{1.00,0.00,1.00}
\definecolor{maroon1}{rgb}{1.00,0.20,0.70}
\definecolor{maroon2}{rgb}{0.93,0.19,0.65}
\definecolor{maroon3}{rgb}{0.80,0.16,0.56}
\definecolor{maroon4}{rgb}{0.55,0.11,0.38}
\definecolor{maroon}{rgb}{0.69,0.19,0.38}
\definecolor{mediumaquamarine}{rgb}{0.40,0.80,0.67}
\definecolor{mediumblue}{rgb}{0.00,0.00,0.80}
\definecolor{mediumorchid}{rgb}{0.73,0.33,0.83}
\definecolor{mediumpurple}{rgb}{0.58,0.44,0.86}
\definecolor{mediumsea}{rgb}{0.24,0.70,0.44}
\definecolor{mediumslate}{rgb}{0.48,0.41,0.93}
\definecolor{mediumspring}{rgb}{0.00,0.98,0.60}
\definecolor{mediumturquoise}{rgb}{0.28,0.82,0.80}
\definecolor{mediumviolet}{rgb}{0.78,0.08,0.52}
\definecolor{midnightblue}{rgb}{0.10,0.10,0.44}
\definecolor{mintcream}{rgb}{0.96,1.00,0.98}
\definecolor{mistyrose}{rgb}{1.00,0.89,0.88}
\definecolor{moccasin}{rgb}{1.00,0.89,0.71}
\definecolor{navajowhite}{rgb}{1.00,0.87,0.68}
\definecolor{navyblue}{rgb}{0.00,0.00,0.50}
\definecolor{navy}{rgb}{0.00,0.00,0.50}
\definecolor{oldlace}{rgb}{0.99,0.96,0.90}
\definecolor{olivedrab}{rgb}{0.42,0.56,0.14}
\definecolor{orange1}{rgb}{1.00,0.65,0.00}
\definecolor{orange2}{rgb}{0.93,0.60,0.00}
\definecolor{orange3}{rgb}{0.80,0.52,0.00}
\definecolor{orange4}{rgb}{0.55,0.35,0.00}
\definecolor{orangered}{rgb}{1.00,0.27,0.00}
\definecolor{orange}{rgb}{1.00,0.65,0.00}
\definecolor{orchid1}{rgb}{1.00,0.51,0.98}
\definecolor{orchid2}{rgb}{0.93,0.48,0.91}
\definecolor{orchid3}{rgb}{0.80,0.41,0.79}
\definecolor{orchid4}{rgb}{0.55,0.28,0.54}
\definecolor{orchid}{rgb}{0.85,0.44,0.84}
\definecolor{palegoldenrod}{rgb}{0.93,0.91,0.67}
\definecolor{palegreen}{rgb}{0.60,0.98,0.60}
\definecolor{paleturquoise}{rgb}{0.69,0.93,0.93}
\definecolor{paleviolet}{rgb}{0.86,0.44,0.58}
\definecolor{papayawhip}{rgb}{1.00,0.94,0.84}
\definecolor{peachpuff}{rgb}{1.00,0.85,0.73}
\definecolor{peru}{rgb}{0.80,0.52,0.25}
\definecolor{pink1}{rgb}{1.00,0.71,0.77}
\definecolor{pink2}{rgb}{0.93,0.66,0.72}
\definecolor{pink3}{rgb}{0.80,0.57,0.62}
\definecolor{pink4}{rgb}{0.55,0.39,0.42}
\definecolor{pink}{rgb}{1.00,0.75,0.80}
\definecolor{plum1}{rgb}{1.00,0.73,1.00}
\definecolor{plum2}{rgb}{0.93,0.68,0.93}
\definecolor{plum3}{rgb}{0.80,0.59,0.80}
\definecolor{plum4}{rgb}{0.55,0.40,0.55}
\definecolor{plum}{rgb}{0.87,0.63,0.87}
\definecolor{powderblue}{rgb}{0.69,0.88,0.90}
\definecolor{purple1}{rgb}{0.61,0.19,1.00}
\definecolor{purple2}{rgb}{0.57,0.17,0.93}
\definecolor{purple3}{rgb}{0.49,0.15,0.80}
\definecolor{purple4}{rgb}{0.33,0.10,0.55}
\definecolor{purple}{rgb}{0.63,0.13,0.94}
\definecolor{red1}{rgb}{1.00,0.00,0.00}
\definecolor{red2}{rgb}{0.93,0.00,0.00}
\definecolor{red3}{rgb}{0.80,0.00,0.00}
\definecolor{red4}{rgb}{0.55,0.00,0.00}
\definecolor{red}{rgb}{1.00,0.00,0.00}
\definecolor{rosybrown}{rgb}{0.74,0.56,0.56}
\definecolor{royalblue}{rgb}{0.25,0.41,0.88}
\definecolor{saddlebrown}{rgb}{0.55,0.27,0.07}
\definecolor{salmon1}{rgb}{1.00,0.55,0.41}
\definecolor{salmon2}{rgb}{0.93,0.51,0.38}
\definecolor{salmon3}{rgb}{0.80,0.44,0.33}
\definecolor{salmon4}{rgb}{0.55,0.30,0.22}
\definecolor{salmon}{rgb}{0.98,0.50,0.45}
\definecolor{sandybrown}{rgb}{0.96,0.64,0.38}
\definecolor{seagreen}{rgb}{0.18,0.55,0.34}
\definecolor{seashell1}{rgb}{1.00,0.96,0.93}
\definecolor{seashell2}{rgb}{0.93,0.90,0.87}
\definecolor{seashell3}{rgb}{0.80,0.77,0.75}
\definecolor{seashell4}{rgb}{0.55,0.53,0.51}
\definecolor{seashell}{rgb}{1.00,0.96,0.93}
\definecolor{sienna1}{rgb}{1.00,0.51,0.28}
\definecolor{sienna2}{rgb}{0.93,0.47,0.26}
\definecolor{sienna3}{rgb}{0.80,0.41,0.22}
\definecolor{sienna4}{rgb}{0.55,0.28,0.15}
\definecolor{sienna}{rgb}{0.63,0.32,0.18}
\definecolor{skyblue}{rgb}{0.53,0.81,0.92}
\definecolor{slateblue}{rgb}{0.42,0.35,0.80}
\definecolor{slategray}{rgb}{0.44,0.50,0.56}
\definecolor{slategrey}{rgb}{0.44,0.50,0.56}
\definecolor{snow1}{rgb}{1.00,0.98,0.98}
\definecolor{snow2}{rgb}{0.93,0.91,0.91}
\definecolor{snow3}{rgb}{0.80,0.79,0.79}
\definecolor{snow4}{rgb}{0.55,0.54,0.54}
\definecolor{snow}{rgb}{1.00,0.98,0.98}
\definecolor{springgreen}{rgb}{0.00,1.00,0.50}
\definecolor{steelblue}{rgb}{0.27,0.51,0.71}
\definecolor{tan1}{rgb}{1.00,0.65,0.31}
\definecolor{tan2}{rgb}{0.93,0.60,0.29}
\definecolor{tan3}{rgb}{0.80,0.52,0.25}
\definecolor{tan4}{rgb}{0.55,0.35,0.17}
\definecolor{tan}{rgb}{0.82,0.71,0.55}
\definecolor{thistle1}{rgb}{1.00,0.88,1.00}
\definecolor{thistle2}{rgb}{0.93,0.82,0.93}
\definecolor{thistle3}{rgb}{0.80,0.71,0.80}
\definecolor{thistle4}{rgb}{0.55,0.48,0.55}
\definecolor{thistle}{rgb}{0.85,0.75,0.85}
\definecolor{tomato1}{rgb}{1.00,0.39,0.28}
\definecolor{tomato2}{rgb}{0.93,0.36,0.26}
\definecolor{tomato3}{rgb}{0.80,0.31,0.22}
\definecolor{tomato4}{rgb}{0.55,0.21,0.15}
\definecolor{tomato}{rgb}{1.00,0.39,0.28}
\definecolor{turquoise1}{rgb}{0.00,0.96,1.00}
\definecolor{turquoise2}{rgb}{0.00,0.90,0.93}
\definecolor{turquoise3}{rgb}{0.00,0.77,0.80}
\definecolor{turquoise4}{rgb}{0.00,0.53,0.55}
\definecolor{turquoise}{rgb}{0.25,0.88,0.82}
\definecolor{violetred}{rgb}{0.82,0.13,0.56}
\definecolor{violet}{rgb}{0.93,0.51,0.93}
\definecolor{wheat1}{rgb}{1.00,0.91,0.73}
\definecolor{wheat2}{rgb}{0.93,0.85,0.68}
\definecolor{wheat3}{rgb}{0.80,0.73,0.59}
\definecolor{wheat4}{rgb}{0.55,0.49,0.40}
\definecolor{wheat}{rgb}{0.96,0.87,0.70}
\definecolor{whitesmoke}{rgb}{0.96,0.96,0.96}
\definecolor{white}{rgb}{1.00,1.00,1.00}
\definecolor{yellow1}{rgb}{1.00,1.00,0.00}
\definecolor{yellow2}{rgb}{0.93,0.93,0.00}
\definecolor{yellow3}{rgb}{0.80,0.80,0.00}
\definecolor{yellow4}{rgb}{0.55,0.55,0.00}
\definecolor{yellowgreen}{rgb}{0.60,0.80,0.20}
\definecolor{yellow}{rgb}{1.00,1.00,0.00}

%% file: Macros_Master.tex
\newcommand{\C}{\mathbb{C}}

\newcommand{\Q}{\mathbb{Q}}
\newcommand{\R}{\mathbb{R}}

\newcommand{\Z}{\mathbb{Z}}
\newcommand{\del}{\partial}

\newcommand{\into}{\hookrightarrow}

\newcommand{\trans}{\pitchfork}

\newcommand{\proj}[2]{\mathbb{#1 P}^{#2}}

\newcommand{\lan}{\left\langle}
\newcommand{\ran}{\right\rangle}
\newcommand{\bgd}{\begin{displaymath}}
\newcommand{\edd}{\end{displaymath}}
\newcommand{\bge}{\begin{equation}}
\newcommand{\ede}{\end{equation}}
\newcommand{\bgea}{\begin{eqnarray*}}
\newcommand{\edea}{\end{eqnarray*}}
\newcommand{\bgeA}{\begin{eqnarray}}
\newcommand{\edeA}{\end{eqnarray}}
\newcommand{\bgc}{\begin{center}}
\newcommand{\edc}{\end{center}}
\newcommand{\ben}{\begin{enumerate}}
\newcommand{\een}{\end{enumerate}}
\newcommand{\bgi}{\begin{itemize}}
\newcommand{\edi}{\end{itemize}}
\newcommand{\hf}{\hspace*{0.5cm}}

\newcommand{\ep}{\varepsilon}
 
\newcommand{\lp}{\bullet}